\documentclass[12pt]{article}

\usepackage{amssymb,amsmath,latexsym}
\oddsidemargin=\evensidemargin \advance\topmargin by-\headheight

\DeclareMathOperator{\aut}{Aut}

\DeclareMathOperator{\ann}{ann}

\DeclareMathOperator{\cyc}{Cyc}

\DeclareMathOperator{\GF}{GF}

\DeclareMathOperator{\GR}{GR}

\DeclareMathOperator{\id}{id}
\def\LI{I_{\rm L}}
\DeclareMathOperator{\im}{im}

\DeclareMathOperator{\md}{mod}
\DeclareMathOperator{\orb}{Orb}

\DeclareMathOperator{\rad}{rad}
\DeclareMathOperator{\rk}{rk}

\DeclareMathOperator{\Span}{span}
\DeclareMathOperator{\supp}{Supp}

\def\UI{I_{\rm U}}

\def\C{{\mathbb C}}

\def\Q{{\mathbb Q}}

\def\Z{{\mathbb Z}}

\def\A{{\cal A}}

\def\H{{\cal H}}
\def\I{{\cal I}}

\def\P{{\cal P}}

\def\S{{\cal S}}
\def\T{{\cal T}}
\def\U{{\cal U}}

\def\wh{\widehat}
\def\wt{\widetilde}

\def\proof{{\bf Proof}.\ }
\def\bull{\vrule height .9ex width .8ex depth -.1ex }

\makeatletter
\renewcommand{\subsection}{\@startsection{subsection}{2}{0mm}{-2mm}{-2mm}
{\bf\normalsize}}

\def\sbsnt#1{\subsection{\hspace{-3mm}#1}}
\makeatother

\newtheorem{formula}{}[section]
\newtheorem{proposition}[formula]{Proposition}
\newtheorem{definition}[formula]{Definition}
\newtheorem{corollary}[formula]{Corollary}
\newtheorem{remark}[formula]{Remark}
\newtheorem{lemma}[formula]{Lemma}
\newtheorem{theorem}[formula]{Theorem}

\newtheorem{example}[formula]{Example}

\def\thrm{\begin{theorem}}
\def\thrml#1{\begin{theorem}\label{#1}}
\def\ethrm{\end{theorem}}
\def\prpstn{\begin{proposition}}
\def\prpstnl#1{\begin{proposition}\label{#1}}
\def\eprpstn{\end{proposition}}
\def\rmrk{\begin{remark}}
\def\rmrkl#1{\begin{remark}\label{#1}}
\def\ermrk{\end{remark}}
\def\dfntn{\begin{definition}}
\def\dfntnl#1{\begin{definition}\label{#1}}
\def\edfntn{\end{definition}}
\def\nmrt{\begin{enumerate}}
\def\enmrt{\end{enumerate}}
\def\tm#1{\item[{\rm (#1)}]}
\def\qtn{\begin{equation}}
\def\qtnl#1{\begin{equation}\label{#1}}
\def\eqtn{\end{equation}}
\def\lmm{\begin{lemma}}
\def\lmml#1{\begin{lemma}\label{#1}}
\def\elmm{\end{lemma}}
\def\crllr{\begin{corollary}}
\def\crllrl#1{\begin{corollary}\label{#1}}
\def\ecrllr{\end{corollary}}
\def\hpthss{\begin{hypothesis}}
\def\hpthssl#1{\begin{hypothesis}\label{#1}}
\def\ehpthss{\end{hypotxesis}}
\def\xmpl{\begin{example}}
\def\xmpll#1{\begin{example}\label{#1}}
\def\exmpl{\end{example}}
\def\css{\begin{cases}}
\def\ecss{\end{cases}}

\begin{document}

\title{Schur rings over a product of Galois rings}
\author{
Sergei Evdokimov \\[-2pt]
\small Steklov Institute of Mathematics\\[-4pt]
\small at St. Petersburg \\[-4pt]
{\tt \small evdokim@pdmi.ras.ru }
\thanks{The work was partially supported by RFFI Grant 07-01-00485.}
\and
Ilia Ponomarenko\\[-2pt]
\small Steklov Institute of Mathematics\\[-4pt]
\small at St. Petersburg \\[-4pt]
{\tt \small inp@pdmi.ras.ru}
\thanks{The work was partially supported by RFFI Grants 07-01-00485, 08-01-00379
and 08-01-00640.}
}
\date{}

\maketitle

\begin{abstract}
The recently developed theory of Schur rings over a finite cyclic group is
generalized to Schur rings over a ring $R$ being a product of Galois rings of
coprime characteristics. It is proved that if the characteristic of $R$ is odd,
then as in the cyclic group case any pure Schur ring over~$R$ is the tensor
product of a pure cyclotomic ring and Schur rings of rank~$2$ over non-fields.
Moreover, it is shown that in contrast to the cyclic group case there are
non-pure Schur rings over $R$ that are not generalized wreath products.
\end{abstract}

\section{Introduction}

In papers \cite{LM96,LM1} K.~H.~Leung and S.~H.~Man proved that any Schur ring
(S-ring) over a finite cyclic group can be constructed from special S-rings
by means of two operations: tensor product and wedge product (as for a background of
S-rings see Section~\ref{011209b}). This theorem supplemented with the normality
theory from~\cite{EP01be} enabled to get a series of strong results in
algebraic combinatorics~\cite{EP01be,EP03be,K04,M01}.\medskip

To generalize the Leung-Man theorem in some way to S-rings over an arbitrary
abelian group, the notion of S-ring over
a commutative ring $R$\,\footnote{If the opposite is not indicated, all rings are supposed to be finite rings
with identity.} was introduced in~\cite{EP03c}; by definition any such ring is
an S-ring over the additive group $R^+$ of the ring $R$ that is invariant
with respect to its multiplicative group $R^\times$. It should be noted that
by the Schur theorem on multipliers any S-ring over a cyclic group of order $n$
can be treated as an S-ring over the ring $R=\Z_n$ of integers modulo~$n$.
We observe that in this case $R$ is the direct product of Galois rings of coprime
characteristics with prime residue fields. Thus it is natural to try to generalize the
Leung-Man theorem to S-rings over the products of arbitrary Galois rings of
coprime characteristics. In this paper such rings are called {\it CG-rings}.\medskip

The first step to constructing the theory of S-rings over rings is to characterize all primitive
S-rings over a ring, i.e. those that have no proper quotients. The
situation is controlled by a generalization of the Burnside-Schur theorem
proved in~\cite{EP03c} (see also Theorem~\ref{201207a}). It turned out that
any primitive S-ring over a ring~$R$ from a quite general class including all CG-rings  is either of rank~$2$, or a cyclotomic ring over $R$ (in the latter case
$R$ is a field). It should be mentioned that the special S-rings in the Leung-Man
theorem belong exactly to one of these two classes. In~\cite{EP09}
the Burnside-Schur theorem was applied to find a complete generalization
of the Leung-Man theorem to S-rings over a Galois ring of odd characteristic.\medskip

On the second step, following the logic of the cyclic group case
we should find a condition for an S-ring over a CG-ring~$R$ to be
the generalized wreath product of S-rings over smaller CG-rings.
(In the cyclic group case this operation  was introduced in~\cite{EP01ae} and produces
exactly the S-rings which are wedge products mentioned above.)
In the case when $R$ is a Galois ring of odd characteristic as well as in
the cyclic case the required condition reduces to the non-purity of the S-ring in question. The
latter means that any of its basic sets intersecting~$R^\times$ is a union of cosets
modulo a fixed non-zero ideal of~$R$ (see Subsection~\ref{021209b}). However,
the following theorem which we prove in Section~\ref{191109a} shows that the
case of an arbitrary CG-ring is more complicated.

\thrml{210809a}
Let $p$ and $q$ be distinct primes, and $d$ and $e$ positive integers such that
$p$ divides $q^e-1$ and $q$ divides $p^d-1$. Set $R=R_p\times R_q$ where
$R_p=\GR(p^2,d)$ and $R_q=\GR(q^2,e)$. Then there exists a non-pure dense\footnote{As
for the definition of density see the beginning of Section~\ref{021209i}.}
S-ring over the CG-ring $R$ that is not a non-trivial generalized wreath product.
\ethrm

The S-rings from Theorem~\ref{210809a} never exist when $R=\Z_n$ because
in this case $e=d=1$ and the hypothesis is obviously not satisfied. It should be
also mentioned that probably all these S-rings are not schurian (as to the
concept of schurity we refer to~\cite{EP09a}). For example, this is the
case when  $(p,d)=(2,2)$ and $(q,e)=(3,1)$.\medskip

In spite of the fact that a complete analog for the non-pure part of the cyclic case theory
can not be reconstructed for general CG-rings (Theorem~\ref{210809a}), some information on
non-pure S-rings over a CG-ring can be obtained. Namely, Theorem~\ref{140409f} in particular
shows that any such ring which is not a non-trivial tensor or generalized wreath
product, must "contain"\ all maximal and minimal ideals of~$R$.
(This theorem is also used to prove the density of pure S-rings.)
The proof of this result occupies Sections~\ref{241109a} and~\ref{021209u}.
As an immediate consequence of Theorem~\ref{231109a} proved there we obtain the
following classification of rational S-rings over a CG-ring, i.e. those
any basic set of which is $R^\times$-invariant.

\thrml{201109a}
Any rational S-ring over a CG-ring is either a non-trivial generalized wreath
product, or a tensor product, one factor of which is an S-ring of rank~$2$.
\ethrm

At the final step we have to characterize pure S-rings over a CG-ring~$R$.
Again as in the cyclic group case we could expect that any such ring is
the tensor product of a pure cyclotomic ring and S-rings of rank~$2$.
However, this is not generally true. For example,
if the characteristic of $R$ is even, then there exist pure dense S-rings  which
are not cyclotomic.~\footnote{Some examples
of such S-rings were found by the authors and will be published elsewhere.}
The following statement shows that this obstacle is a unique one.

\thrml{270409a}
Any pure S-ring over a CG-ring of odd characteristic is the tensor product of a
pure cyclotomic ring and S-rings of rank~$2$ over non-fields.
\ethrm

Theorem~\ref{270409a} is an immediate consequence of Theorems~\ref{160409a}
and~\ref{130209b} proved in Sections~\ref{021209i} and~\ref{271109a} respectively.
In the proofs of these theorems we use the duality theory for S-rings over
a CG-ring developed in Subsections~\ref{021109k} and~\ref{021109l}. In the first
case (Theorem~\ref{160409a}) this theory shows that the property of an S-ring to be
a non-trivial generalized wreath or tensor product is preserved
under duality, 
and enables us to interchange minimal and maximal ideals of the
ring. This reduces the proof to the dense case. In the second case
(Theorem~\ref{130209b}) we use another fact from this theory: an S-ring and its dual
are cyclotomic or not simultaneously. This enables us to prove that any
dense pure S-ring over a CG-ring contains a pure cyclotomic S-ring. The rest
of the proof is heavily based on Theorem~\ref{250303b} which applies to separate special
pure sets by means of characters of the additive group of the underlying ring.\medskip

It should be stressed that this work was essentially influenced by papers~\cite{LM96} and~\cite{M94}.
However, there is a great difference between the cyclic group case and the general CG-ring case.
Namely, in the former case the projection of any pure subgroup of the multiplicative group of
the ring on at least one local component is also pure. On the contrary, in the latter one this is
not true, e.g. for the CG-rings satisfying the hypothesis of Theorem~\ref{210809a}.
In fact, this is the only reason why the non-pure part of the theory can not be done
properly. On the other hand, it is quite surprising that the pure part is
completed by Theorem~\ref{201109a}.\medskip

Concerning finite rings and permutation groups we refer to~\cite{MD74} and~\cite{DM}.
For the reader convenience we collect the basic facts on S-rings over abelian
groups, on CG-rings and on S-rings over them in Sections~\ref{011209b},
\ref{021109d} and \ref{021209x} respectively. In Section~\ref{021209w} we study
pure sets in CG-rings. One of the main results here is the separation theorem
(Theorem~\ref{250303b}); the proof of it is given in Section~\ref{271109b}.\medskip

{\bf Notation.}
As usual by $\Z,\Q,\C$ we denote the ring of integers, the ring of rationals
and the field of complex numbers respectively.

For a prime $p$ the $p$th part of a positive integer $n$ is denoted by $n_p$.

For a commutative ring~$R$ with identity we denote by $R^\times$ and $\rad(R)$ the
multiplicative group of~$R$ and the radical of~$R$ respectively.

The set of all (resp. all maximal, all minimal) ideals of $R$ is denoted by
$\I(R)$ (resp. $\I_{max}(R)$, $\I_{min}(R)$).

Given $I\in\I(R)$ we denote by $I^+$ the additive group of~$I$, and by $\pi_I$
the natural epimorphism from~$R$ to~$R/I$.

For a set $X\subset R$ we denote by $\UI(X)$ the smallest ideal of $R$
containing~$X$ and by $\LI(X)$ the largest ideal $I$ of $R$ such that
$X+I=X$ or, equivalently, that $X$ is a union of $I$-cosets. Also we set
$$\ann(X)=\{r\in R:\ rX=\{0\}\}$$
and write $\ann(r)$ instead of $\ann(\{r\})$ for $r\in R$.

Let $G=\prod_{p\in\P}G_p$ be a finite abelian group where $\P=\P(G)$ is the set
of all primes dividing $|G|$ and $G_p$ is the Sylow $p$-subgroup of~$G$.
For $Q\subset\P$ the $Q$-projection of $x\in G$ (resp. $X\subset G$) is denoted
by $x_Q$ (resp. $X_Q$). When $Q=\{p\}$ we omit the braces and we write
$Q'$ instead of $\P\setminus Q$. For an arbitrary set $Q$ of primes we
put $x_Q=x_{Q\cap\P}$ and $X_Q=X_{Q\cap\P}$.

For a subset $X$ of a group $G$ we set $X^\#=X\setminus\{1_G\}$.

The group ring of a group $G$ over an arbitrary ring $R$ is denoted by $RG$.
The element of~$RG$ that is equal to the sum of all elements of a set $X\subset G$ is
denoted by $\xi(X)$. 
The support of $\xi\in RG$ is denoted by $\supp(\xi)$. 
The componentwise multiplication in $RG$ is denoted by~$\circ$.

\section{S-rings over groups}\label{011209b}
\sbsnt{Definition and properties.}\label{081209i}
Let $G$ be a finite group. A subring~$\A$ of the group ring~$\Z G$ is
called a {\it Schur ring} ({\it S-ring}, for short) over~$G$ if it has a
(uniquely determined) $\Z$-basis consisting of the elements
$\xi(X)=\sum_{x\in X}x$ where $X$ runs over a family $\S=\S(\A)$ of pairwise
disjoint non-empty subsets of~$G$ such that
$$
\{1\}\in\S,\quad
\bigcup_{X\in\S}X=G\quad
\textstyle{\rm and}\quad
X\in\S\ \Rightarrow\ X^{-1}\in\S.
$$
We call the elements of $\S$ the {\it basic} sets of~$\A$ and denote by $\S^*(\A)$
the set of all unions of them and by~$\H(\A)$ the set of all subgroups of~$G$
in $\S^*(\A)$. The elements of $\S^*(\A)$ and $\H(\A)$ are called
{\it $\A$-subsets of~$G$} (or {\it $\A$-sets}) and {\it $\A$-subgroups of~$G$}
respectively. It is easily seen that $XY$ is an $\A$-set whenever so are~$X$
and~$Y$. For an $\A$-set $X$ we put
$$
\S_X=\{X'\in\S:\ X'\subset X\}.
$$
The number $\rk(\A)=\dim_\Z(\A)$ is called the {\it rank} of~$\A$. If
$\S^*(\A)\subset\S^*(\A')$ where $\A'$ is an S-ring over~$G$, then
we write $\A\le\A'$.\medskip

\lmml{090608a}
Let $\A$ be an S-ring over a group $G$, $H\in\H(\A)$ and $X\in\S(\A)$. Then
the cardinality of the set $X_{H,x}=X\cap Hx$ does not depend on $x\in X$.
\elmm
\proof See \cite[p.21]{EP09}.\bull\medskip

Let $H\in\H(\A)$. Then $\S_H$ where $\S=\S(\A)$, is the set of basic sets of an
S-ring over the group $H$. This S-ring is denoted by~$\A_H$. If the
group $H$ is normal and $\pi:G\to G/H$ is the quotient epimorphism, then
$\S_{G/H}=\pi(\S)$ is the set of basic sets of an S-ring over the group $G/H$.
This S-ring is denoted by~$\A_{G/H}$.\medskip

If $\A_1$ and $\A_2$ are S-rings over groups $G_1$ and $G_2$ respectively, then
the subring $\A=\A_1\otimes \A_2$ of the ring $\Z G_1\otimes\Z G_2=\Z G$ where
$G=G_1\times G_2$, is an S-ring over the group $G$ with
$$
\S(\A)=\{X_1\times X_2: X_1\in\S(\A_1),\ X_2\in\S(\A_2)\}.
$$
It is called the {\it tensor product} of $\A_1$ and $\A_2$.

\lmml{130209d}
Let $G_1$ and $G_2$ be groups, and $\A$ an S-ring over the group $G=G_1\times G_2$.
Suppose that $G_1,G_2\in\H(\A)$. Then $\pi_i(X)\in\S(\A)$ for all $X\in\S(\A)$
where $\pi_i$ is the projection of $G$ on $G_i$, $i=1,2$. In particular,
$\A\ge\A_{G_1}\otimes \A_{G_2}$.
\elmm
\proof Let $X\in\S(\A)$ and $i\in\{1,2\}$. Then obviously $G_{3-i}X\cap G_i=\pi_i(X)$. By
the hypothesis this implies that $\pi_i(X)\in\S^*(\A)$. To complete the proof it
suffices to note that if $\pi_i(X)$ is the disjoint union of non-empty $\A$-sets
$Y$ and $Z$, then $X$ is the disjoint union of the non-empty $\A$-sets
$X\cap G_{3-i}Y$ and $X\cap G_{3-i}Z$ which is impossible because
$X\in\S(\A)$.\bull\medskip

Let $\A$ be an S-ring over a group $G$ and let $L,U$ be $\A$-subgroups of $G$
such that $L\le U$ and $L$ is normal in~$G$. Following~\cite{EP01be} we say that
$\A$ satisfies the {\it $U/L$-condition} if
$$
LX=XL=X,\qquad X\in\S(\A)_{G\setminus U}.
$$
If, moreover, $L\ne\{1\}$ and $U\ne G$, we say that $\A$ satisfies the
$U/L$-condition {\it non-trivially}.\medskip

An S-ring~$\A$ satisfying the $U/L$-condition was called in~\cite{LM96,LM1} the
wedge product of the S-rings $\A_U$ and $\A_{G/L}$. It should be noted that
the authors in~\cite{EP01ae} independently introduced the external operation of the
generalized wreath product of two S-rings which produces exactly the S-rings
satisfying the $U/L$-condition.\medskip

The following important theorem goes back to I.~Schur and H.~Wielandt (see \cite[Ch.~IV]{W64});
as to the formulation given here we refer to~\cite{EP03c}. Below for
$X\subset G$, $m\in\Z$, and a prime $p$ we set
$$
X^{(m)}=\{x^m:\ x\in X\},\quad
X^{[p]}=\{x^p:\ x\in X,\ |xH\cap X|\not\equiv 0\pmod p\}
$$
where $H=\{g\in G:\ g^p=1\}$.

\thrml{261009b}
Let $G$ be a finite abelian group and $\A$ an S-ring over $G$. Then for any $X\in\S(\A)$ the
following statements hold:
\nmrt
\tm{1} $X^{(m)}\in\S(\A)$ for any integer $m$ coprime to $|G|$,
\tm{2} $X^{[p]}\in\S^*(\A)$ for any prime $p$ dividing $|G|$.\bull
\enmrt
\ethrm

\sbsnt{Duality.}\label{151209u}
Let $\A$ be an S-ring over a finite abelian group~$G$ and $\wh G$ the
group dual to~$G$, i.e. the group of all irreducible $\C$-characters of~$G$.
Given $S\subset G$ and $\chi\in\wh G$ set
\qtnl{201108b}
\chi(S)=\sum_{s\in S}\chi(s).
\eqtn
Characters $\chi_1,\chi_2\in\wh G$ are called equivalent if
$\chi_1(S)=\chi_2(S)$ for all $S\in\S(\A)$. Denote by $\wh\S$ the set of
classes of this equivalence relation.
Then the submodule of $\Z \wh G$ spanned by the elements $\xi(X)$, $X\in\wh\S$,
is an S-ring over~$\wh G$ (see~\cite[Theorem~6.3]{BI}). This ring is called
{\it dual} to~$\A$ and is denoted by~$\wh\A$. Obviously, $\S(\wh\A)=\wh\S$.
Moreover, $\rk(\wh\A)=\rk(\A)$ and
\qtnl{311008a}
\H(\wh\A)=\{H^\bot:\ H\in \H(\A)\}
\eqtn
where $H^\bot=\{\chi\in\wh G:\ H\le\ker(\chi)\}$. It is also true that the S-ring
dual to $\wh\A$ is equal to~$\A$. The following theorem was proved in~\cite{EP09}.

\thrml{260808a}
Let $\A$ be an S-ring over an abelian group $G$. 
Then $\A$ satisfies the $U/L$-condition if and only if $\wh\A$ satisfies
the $L^\bot/U^\bot$-condition.\bull
\ethrm

Some more properties of the dual S-ring are contained in the following statement.

\thrml{140509a}
Let $\A$ be an S-ring over an abelian group $G$. Then
\nmrt
\tm{1} $\wh{\A_H}=\wh\A_{\wh G/H^\bot}$ and $\wh{\A_{G/H}}=\wh\A_{H^\bot}$
for any $H\in\H(\A)$,
\tm{2} $\A=\A_1\otimes\A_2$ if and only if $\wh\A=\wh\A_1\otimes\wh\A_2$.
\enmrt
\ethrm
\proof To prove statement~(1) it suffices to verify the second equality. It is
easily seen that $\chi(X)=a\chi(XH)$ for all $\chi\in H^\bot$ and $X\in\S(\A)$
where $a$ is a positive rational. So for any $\chi_1,\chi_2\in H^\bot$ we have
$$
\chi_1(XH)=\chi_2(XH)\ \Leftrightarrow\ \chi_1(X)=\chi_2(X)
$$
and we are done by the definition of the dual S-ring. To prove statement~(2)
suppose that $\A=\A_1\otimes\A_2$. Then $G=G_1\times G_2$ where $G_1,G_2\in\H(\A)$,
and $\A_1=\A_{G_1}=\A_{G/G_2}$ and $\A_2=\A_{G_2}=\A_{G/G_1}$ (see Lemma~\ref{130209d}).
Therefore from statement~(1) it follows that
$$
\wh\A_i=\wh{\A_{G_i}}=\wh{\A_{G/G_{3-i}}}=\wh\A_{G_{3-i}^\bot},\qquad i=1,2
$$
(we identify $\wh{G_i}$ and $G_{3-i}^\bot$). Since obviously $\wh G=G_1^\bot\times G_2^\bot$
and $G_1^\bot,G_2^\bot\in\H(\wh{\A})$ (see~(\ref{311008a})), this implies that
$\wh\A\ge \wh\A_1\otimes\wh\A_2$. Since also
$\rk(\wh{\A_i})=\rk(\A_i)$, $i=1,2$, and $\rk(\A)=\rk(\wh\A)$, we conclude
that $\wh\A=\wh\A_1\otimes\wh\A_2$.\bull

\section{CG-rings}\label{021109d}

Throughout the rest of the paper under a ring we mean a finite commutative ring with
identity.

\sbsnt{Products of Galois rings.}\label{021109dd}
Following \cite[Section XVI]{MD74} 
a local ring $R$ is called {\it Galois} if it is a Galois extension of the prime
ring $\Z_{p^n}$ for some prime~$p$ and positive integer~$n$, or equivalently if
$\rad(R)=pR$.
Given positive integers $n,d$ there exists a unique (up to isomorphism) Galois
ring of characteristic $p^n$ with the residue field of order $q=p^d$; it is
denoted by $\GR(p^n,d)$. We observe that
$$
\GR(p,d)\cong\GF(p^d),\qquad \GR(p^n,1)\cong\Z_{p^n}.
$$
Each ideal of the  Galois ring $\GR(p^n,d)=R$ other than $R$ is of the
form $p^iR$, $i=1,\ldots,n$, and the corresponding quotient ring is isomorphic
to $\GR(p^i,d)$. It is known that $R^+$ is a homocyclic $p$-group of rank $d$ and
exponent $p^n$, i.e. it is isomorphic to the direct product of $d$ cyclic groups
of order $p^n$. Moreover,
\qtnl{250706b}
R^\times=\T\times\U
\eqtn
where $\T$ is the Teichm\"uller group and $\U$ is the group of principal units.
The groups $\T$ and $\U=1+\rad(R)$ are a cyclic group of order $q-1$ and an
abelian $p$-group respectively. If $p$ is odd, then the group $\U$ is homocyclic
of rank~$d$ and exponent $p^{n-1}$.\medskip

Let $R$ be a ring and $\P=\P(R)$. It is well known (see e.g. \cite[Theorem~6.2]{MD74})
that there is a decomposition
$$
R=\prod_{p\in\P}R_p
$$
where $R_p$ is the {\it $p$-component} of~$R$, i.e. the subring of~$R$ such that
$(R_p)^+$ is the Sylow $p$-subgroup of~$R^+$. Moreover, each $R_p$ is the direct
product of local rings the characteristic of each of which is a power of~$p$.
For any $Q\subset\P$ the set $R_Q$ (defined in Notation) equals the product of all
rings $R_p$ with $p\in Q$.

\dfntn
We say that $R$ is a CG-ring (componentwise Galois ring) if $R_p$ is a Galois ring
for all $p\in\P$.
\edfntn

Obviously, the ring with one element as well as any Galois ring is CG. The
characteristic $c$ of a CG-ring $R$ equals the product of the characteristics
$c_p$ of its components $R_p$, $p\in\P$. It is easily seen that
\qtnl{161109c}
\I(R)=\{mR:\ m\ \text{divides}\ c\}.
\eqtn
In particular, the minimal and maximal ideals of $R$ are exactly those $mR$
for which respectively $m=c/p$ and $m=p$ where $p\in\P$. Throughout the paper we
denote by $I_0$ the sum of minimal ideals in the components and set $I_{0,p}=(I_0)_p$.
Clearly,
\qtnl{151209a}
I_{0,p}=(c_p/p)R_p,\qquad p\in\P.
\eqtn
One can see that given $I\in\I(R)$ the set $1+I$ is a subgroup of $R^\times$ if
and only if $I_p\ne R_p$ for all $p\in\P(I)$.\medskip

The class of all CG-rings is closed with respect to taking quotients. Moreover,
the following equality holds:
\qtnl{161109b}
\pi_I(R)^\times=\pi_I(R^\times),\qquad I\in\I(R).
\eqtn
Let $a\in R$. Then the mapping $x\mapsto ax$, $x\in R$,
induces an $R$-module epimorphism from $R$ onto $I=aR$ the kernel of which
coincides with $\ann(a)$. Therefore the ring $R/\ann(a)$ and the ideal $I$ are
isomorphic as $R$-modules. This enables us to define a ring structure on
the ideal~$I$. The corresponding ring $R_{I,a}$ has $a$ as identity and
\qtnl{161109a}
f_{I,a}:R\to R_{I,a},\quad x\mapsto ax
\eqtn
is a ring epimorphism with the kernel $\ann(a)$. It should be noted that for any $u\in R^\times$
the rings $R_{I,a}$ and $R_{I,ua}$ are isomorphic.
When $a=m\cdot 1$ where $m$ is a positive integer dividing
the characteristic of~$R$, we set $R_I=R_{I,a}$
and $f_I=f_{I,a}$. From~(\ref{161109c}), (\ref{161109b})
and~(\ref{161109a}) it follows that in this case
\qtnl{261109u}
(R_I)^\times=mR^\times.
\eqtn

In this paper we consider the permutation group induced by the action of the
group $R^\times$ on the set $R$ by multiplication. It is a subgroup of the
group $\aut(R^+)$ that leaves any ideal of $R$ fixed. Moreover, the orbits of
this group are regular and $R^\times$ is the only faithful one.

\sbsnt{Duality.}\label{021109k}
Let $R$ be a CG-ring of characteristic $c$ and $\P=\P(R)$.
For each $p\in\P$ denote by $\chi_p$ a character~ of the group $(R_p)^+$ such that
$\im(\chi_p)$ contains a primitive $c_p$th root of unity, and by $\wh R_p$
the Galois ring {\it dual} to~$R_p$ with respect to $\chi_p$  (see~\cite{EP09}).
The CG-ring
$$
\wh R=\prod_{p\in\P}\wh R_p
$$
is called {\it dual} to~$R$ with respect to the character~$\chi=\prod_p\chi_p$.
Clearly, $\wh R^+=\wh{R^+}$ is the group dual to the group~$R^+$, and
\qtnl{010908a}
(mR)^\bot=(c/m){\wh R}
\eqtn
where $m$ is a divisor of $c$. Moreover, $\wh R=\{\chi^{(r)}:\ r\in R\}$ where $\chi^{(r)}$ is
the character of $R^+$ such that $\chi^{(r)}(x)=\chi(rx)$, $x\in R$, and the
multiplication in $\wh R$ is defined by the formula $\chi^{(r)}\chi^{(s)}=\chi^{(rs)}$,
$r,s\in R$. Thus $\chi$ is the identity of this ring and
\qtnl{301008b}
\wh R^\times=\{\chi^{(r)}:\ r\in R^\times\}.
\eqtn
The mapping $r\mapsto\chi^{(r)}$ is a ring isomorphism from $R$ onto~$\wh R$;
the image of a group $K\le R^\times$ with respect to this isomorphism is
denoted by~$\wh K$.

\section{Purity in CG-rings}\label{021209w}
Let $R$ be a ring.
Following~\cite{EP07r} a non-empty set $X\subset R$ is called {\it pure}
if $\LI(X)=0$. This means that given $I\in\I(R)$ the equality $X+I=X$ implies
$I=0$. Thus a non-empty set is non-pure if and only if it is a union of $I$-cosets for some
non-zero ideal~$I$.
It is clear that
\qtnl{281008a}
\LI(X)=\LI(uX),\quad u\in R^\times,
\eqtn
and hence the sets $X$ and $uX$ are pure or not simultaneously. Therefore any subset
of a pure orbit of a subgroup of $R^\times$ is also pure. It should be also
noted that the purity of a group $K\le R^\times$, is equivalent to
the fact that $K$ does not contain any subgroup of the form $1+I$ where
$I$ is a non-zero ideal of $R$.\medskip

We note that generally the purity of a set is not preserved under taking quotients
or multiplying by integer. For instance, let
$$
R=\Z_8,\quad X=\{-1,1\},\quad I=4R.
$$
Then the set $X$ is pure whereas the sets $\pi_I(X)$, $2X$ are not.
However, for CG-rings of odd characteristic the situation is controlled as follows.

\thrml{090209a}
Let $R$ be a CG-ring of odd characteristic, $K\le R^\times$ a pure group and $J$
an ideal of $R$. Suppose that $J_p\ne R_p$ for all $p\in\P$.
Then the group $\pi_J(K)$ and the set $c_JK$ are pure where $c_J$ is the
characteristic of~$J$.
\ethrm
\proof Since the ring $R/J$ and the ideal $c_JR$ are isomorphic as $R$-modules,
it suffices to verify only that the set $\pi_J(K)$ is pure. To do this we start
with two observations. First, the condition $J_p\ne R_p$ for all $p\in\P$,
implies that the set $1+J$ is a subgroup of the group $R^\times$. Therefore
there is a canonical isomorphism
\qtnl{121109a}
\pi_J(L)\cong L/(L\cap(1+J)),\qquad L\le R^\times.
\eqtn
Second, denote by $r_p$ the degree of the residue field of the ring~$R_p$
over the prime subfield, $p\in\P$. Then $r_p=\rk(\U_p)$
whenever $R_p$ is not a field where $\U_p$ is the group of principal units
of the ring~$R_p$. We claim that
for any $L\le R^\times$ the following equivalence holds
\qtnl{121109b}
L\ \text{is pure}\quad\Leftrightarrow\quad\rk(L\cap\U_p)<r_p\ \text{for all}\ p\in\P.
\eqtn
To prove (\ref{121109b}) we observe that since $L$ is a group, given an ideal $I\in\I(R)$ we
have $L=L+I$ if and only if $L\ge 1+I$. Therefore the group $L$ is non-pure if and only if
$L\ge 1+I$ for some non-zero ideal $I$ of~$R$ (here $0\ne I_p\ne R_p$ for all $p\in\P(I)$).
The latter is equivalent to the existence of $p\in\P$ such that
$L\ge 1+I_{0,p}$ (here $I_{0,p}\ne R_p$, see above). However, this means that
$$
\rk(1+I_{0,p})=\rk(\U_p)=r_p\quad\text{and}\quad L\cap\U_p\ge 1+I_{0,p},
$$
i.e. $\rk(L\cap\U_p)\ge \rk(1+I_{0,p})=r_p$.\medskip

To complete the proof we note that the rank of an abelian group does not
increase under taking quotients. Therefore from~(\ref{121109a})
we obtain
$$
\rk(\pi_J(K)\cap\pi_J(\U_p))=
\rk((K\cap\U_p)/(K\cap\U_p\cap (1+J))\le
\rk(K\cap\U_p)<r_p
$$
for all $p\in\P$. Since the degree of the residue field of the ring $\pi_J(R_p)$
equals~$r_p$ and the group of principal units of the latter ring coincides with
$\pi_J(\U_p)$, we are done by~(\ref{121109b}).\bull\medskip

The following important statement which will be used in Section~\ref{271109a} is
a kind of separation theorem for pure subsets in CG-rings
of an arbitrary characteristic.

\thrml{250303b}
Let $R$ be a CG-ring, and $S,S'\subset R$ distinct subsets in a pure orbit of a
subgroup of $R^\times$, $S\ne\emptyset$. Then
\nmrt
\tm{1} there exists $\chi\in\wh R^\times$ such that $\chi(S)\ne\chi(S')$,
\tm{2} given $\chi\in\wh R^\times$ there exists $r\in R^\times$ such that
$\chi(rS)\ne 0$.
\enmrt
\ethrm

The proof of this theorem is given in Section~\ref{271109b}.

\section{S-rings over a CG-ring}\label{021209x}
\sbsnt{Definition and properties.}
Let $R$ be a ring and  $\A$ an S-ring over the group~$R^+$. In accordance
with~\cite{EP03c} we say that $\A$ is an {\it S-ring over $R$} if it is invariant
with respect to the action of the group $R^\times$ on $\Z R^+$ by multiplication,
or, equivalently, if for any $u\in R^\times$
\qtnl{171109a}
X\in\S(\A)\quad\Rightarrow\quad uX\in\S(\A).
\eqtn
One can see that this is the case if $\rk(\A)=2$ or $\S(\A)=\orb(K,R)$ for some $K\le R^\times$.
In the latter case $\A$ is called {\it cyclotomic} and denoted by $\cyc(K,R)$.
Clearly,
$$
\cyc(K,R)=\Span\{\xi(X):\ X\in\orb(K,R)\}.
$$
Since the group $K$ acts semiregularly on~$R^\times$, the cyclotomic rings are in
1-1 correspondence to the subgroups of~$R^\times$. The following result is a special
case of~\cite[Corollary~2.3]{EP09}.

\thrml{230810b}
Let $\A$ be an S-ring over a ring $R$. Then any basic set of~$\A$ contained in $R^\times$
is a subgroup of $R^\times$. In particular, any S-ring over a field is a cyclotomic one.\bull
\ethrm

Let $\A$ be an S-ring over the ring $R$. Set
$$
\I(\A)=\I(R,\A)=\{I\in\I(R):\ I\in\S^*(\A)\}.
$$
The elements of $\I(\A)$ are called {\it $\A$-ideals} of~$R$. It was proved in
\cite[Theorem~2.6]{EP09} that if a ring $R$ is generated by the units, then
\qtnl{171109c}
\UI(X),\LI(X)\in\I(\A),\qquad X\in\S^*(\A).
\eqtn
Since the assumption is true for every local ring, the conclusion holds when
$R$ is a CG-ring. From now on we assume that $R$ is a CG-ring.\medskip

Let $\A$ be an S-ring over the ring $R$ and $I\in\I(\A)$. Since the ring $\A$
is $R^\times$-invariant, the S-rings $\A_{R/I}=\A_{R^+/I^+}$ and
$\A_I=\A_{I^+}$ over the groups~$R^+/I^+$ and~$I^+$ are $(R/I)^\times$- and
$(R_I)^\times$-invariant respectively. 
Thus they are S-rings over the rings~$R/I$ and~$R_I$.
\medskip

The S-ring $\A$ is called {\it dense} if $\I(\A)=\I(R)$.
From (\ref{161109c}) it follows that the set-difference between an ideal of $R$
and the union of all proper ideals in it is an orbit of the group~$R^\times$.
Therefore
\qtnl{021109b}
R^\times r\in\S^*(\A),\qquad r\in R.
\eqtn
We say that $\A$ is {\it $R$-primitive} if  $0$ and $R$ are the only
$\A$-ideals of~$R$. The following theorem is a special case of the main
result of~\cite{EP03c}.

\thrml{201207a}
Let $R$ be a CG-ring. Then any $R$-primitive S-ring is either of rank~$2$ or
cyclotomic.\footnote{In the conditions of the theorem "$R$-primitivity"\ is
equivalent to "quasiprimitivity"\ in sense of \cite{EP03c} (see the remark
before Theorem~1.3 of that paper).} In the latter case, $R$ is a field.\bull
\ethrm

The following statement will be used in the proof of Theorem~\ref{011009a}.
(In fact, the proof shows that the conclusion is true under a weaker
assumption that the $p$-component of the ring~$R$ is a Galois ring.)

\thrml{261009a}
Let $R$ be a CG-ring of characteristic $c$ and $\A$ an S-ring over~$R$.
Then for any $p\in\P(R)$ and $X\in\S(\A)$ such that $R_p^\times X=X$, the
following statements hold:
\nmrt
\tm{1} $(X^{[p]})_p=\{0\}$,
\tm{2} $\LI(X)_p\ne 0$ if and only if $X^{[p]}=\emptyset$,
\tm{3} if $\LI(X)_p=0$, then either $X_p=\{0\}$ or $\LI(X\cup Y)_p\ne 0$ with
$Y=(X^{[p]})^{(m)}$ where $m$ is an integer coprime to $c$ such that
$mp\equiv 1\,(\md c_{p'})$.
\enmrt
\ethrm
\proof Let $p\in\P(R)$ and $X\in\S(\A)$ be such that $R_p^\times X=X$. Take $x\in X$.
Then $R_p^\times x\subset X$. On the other hand, set $J=I_{0,p}$ (see~(\ref{151209a})). Since
$R$ is a CG-ring, $R_p$ is a Galois ring, and hence the set $x_pR_p^\times$ is
either a $J$-coset, or $J^\#$, or $\{0\}$ (in the latter two cases the order of
$x_p$ is $p$ and $1$ respectively). Thus
\qtnl{261009c}
X_{x,J}\ne x+J\quad\Leftrightarrow\quad
X_{x,J}=x_{p'}+J^\#\ \text{or}\ X_{x,J}=\{x_{p'}\}.
\eqtn
From the hypothesis of the theorem it follows that the set of all
elements of order $p$ in $R$ coincides with $J^\#$. Since $X$ is a union of all
sets $X_{x,J}$, $x\in X$, we conclude by~(\ref{261009c}) and statement~(2) of Theorem~\ref{261009b} that
$X^{[p]}=pX'$ where $X'$ is the set of all $x\in X$ for which the
condition in the right-hand side of (\ref{261009c}) holds (here we use this
theorem for additively written group $G=R^+$; in this case $H=J$ and
$xH\cap X=X_{J,x}$). This proves statement~(2), and due to the obvious
equality $(pX')_p=\{0\}$ also statement~(1). To prove statement~(3)
suppose that $\LI(X\cup Y)_p=0$. Then since
\qtnl{261009d}
Y=(X^{[p]})^{(m)}=(pX')^{(m)}=mpX'=(X')_{p'},
\eqtn
from (\ref{261009c}) it follows that there exists $x\in X'$ for which
$X_{x,J}=\{x_{p'}\}$. However, in this case $x=x_{p'}$, and hence $x\in X\cap Y$.
Taking into account that $X\in\S(\A)$ and $Y\in\S^*(\A)$ (Theorem~\ref{261009b}), we
conclude that $X\subset Y$. This implies by (\ref{261009d}) that
$X_p\subset Y_p=\{0\}$ which completes the proof.\bull

\sbsnt{Pure S-rings and generalized wreath products.}\label{021209b}
Let $\A$ be an S-ring over a CG-ring $R$. Due to (\ref{281008a}) the ideal
$\LI(X)$ does not depend on the set $X\in\S(\A)$ such that
$X\cap R^\times\ne\emptyset$. We denote this ideal by $\LI(\A)$.

\dfntn
The S-ring $\A$ is {\it pure} if $\LI(\A)=0$.
\edfntn

It follows that $\A$ is pure if and only if some (and hence any) $X\in\S(\A)$,
$X\cap R^\times\ne\emptyset$, is pure. It is easily seen that S-rings of rank~$2$
and cyclotomic rings $\cyc(K,R)$ with pure groups $K\le R^\times$, are pure.

\thrml{021109a}
Let $R$ be a CG-ring of odd characteristic and
$\A$ a pure S-ring over~$R$. Suppose that $\I_{max}(R)\subset\I(\A)$
and $J$ is an $\A$-ideal of $R$ such that $J_p\ne R_p$ for all $p\in\P(R)$.
Then the S-ring $\A_{R/J}$ is also pure.
\ethrm
\proof Let $Y\in\S(\A_{R/J})$ be such that $\pi(1)\in Y$ where $\pi=\pi_J$.
Then it suffices to verify that the set $Y$ is pure. To do this denote by
$X$ the basic set of~$\A$ containing~$1$. Then obviously
$Y=\pi(X)$. Moreover, since by the hypothesis $\I_{max}(R)\subset\I(\A)$, it follows that
$R^\times$ is an $\A$-set, and hence $X\subset R^\times$. So $X$ is a subgroup
of~$R^\times$ by Theorem~\ref{230810b}. However, the set~$X$
is pure because so is the S-ring~$\A$. Therefore the set $Y=\pi(X)$ is pure
by Theorem~\ref{090209a} and we are done.\bull\medskip

In \cite[Theorem~2.8]{EP09} it was proved that an S-ring over a local ring is non-pure
if and only if it is a non-trivial generalized wreath product. However, this is not true for
the S-rings over an arbitrary  CG-ring (see Section~\ref{191109a}). Since most of the results
in this paper are formulated in terms of that product, we recall its definition here.

\dfntnl{deflu}
We say that $\A$ is a generalized wreath product if there exist
$\A$-ideals~$I$ and~$J$ such that
\qtnl{041208b}
J\subset\LI(X)\cap I,\qquad X\in\S(\A)_{R\setminus I}.
\eqtn
In this case we also say that $\A$ satisfies the $I/J$-condition.
\edfntn

If both ideals $I$ and $J$ are proper, we say that the generalized wreath product
is {\it non-trivial}, or that the S-ring~$\A$ satisfies the $I/J$-condition
{\it non-trivially}.\medskip

It should be noted that the S-ring $\A$ satisfying the $I/J$-condition satisfies
also the $I^+/J^+$-condition as defined in Subsection~\ref{081209i}. Therefore
in the sense of~\cite{EP01ae} in this case the S-ring $\A$ is the (standard)
generalized wreath product of the S-rings
$\A_I$ and $\A_{R/J}$ over the groups $I^+$ and $(R/J)^+$ respectively.
Moreover, the latter S-ring can be treated as an S-ring over the ring $R/J$
whereas the former one can be treated as an S-ring over the ring~$R_I$ defined
in Subsection~\ref{021109dd}.\medskip

\sbsnt{Duality.}\label{021109l}
Let $\wh R$ be the ring dual to a CG-ring~$R$ with respect to a character~$\chi$.
We observe that since the set $\wh R^\times$ does not depend on the choice
of~$\chi$, any S-ring over this ring  is also an S-ring over the ring dual
to $R$ with respect to any other character belonging to $\wh R^\times$.

\thrm
Let $\A$ be an S-ring over a CG-ring $R$ and $\wh\A$ the S-ring over
the group~$\wh R^+$ that is dual to~$\A$. Then $\wh\A$ is an S-ring over the
ring~$\wh R$.
\ethrm
\proof Suppose that $\chi^{(s)}$ and $\chi^{(t)}$ belong to the same basic set
of~$\wh\A$ where $s,t\in R$. Then given $r\in R^\times$ we have
$\chi^{(s)}(rS)=\chi^{(t)}(rS)$, or equivalently
$$
\chi^{(rs)}(S)=\chi^{(rt)}(S),\quad S\in\S(\A).
$$
Since $\chi^{(rs)}=\chi^{(r)}\chi^{(s)}$ and $\chi^{(rt)}=\chi^{(r)}\chi^{(t)}$,
this implies that the characters $\chi^{(r)}\chi^{(s)}$ and $\chi^{(r)}\chi^{(t)}$
belong to the same basic set of~$\wh\A$ for all $r\in R^\times$. Thus
the required statement follows from~(\ref{301008b}).\bull\medskip

Let $c$ be the characteristic of the ring $R$ (and hence of the ring $\wh R$). Then
$\I(R)=\{mR:\ m$ divides $c\}$ and $\I(\wh R)=\{m\wh R:\ m$
divides $c\}$. Therefore by~(\ref{010908a}) and~(\ref{311008a}) we have
\qtnl{010908b}
\I({\wh\A})=\{I^\bot:\ I\in\I(\A)\}.
\eqtn
Thus $\A$ is $R$-primitive if and only if $\wh\A$ is $\wh R$-primitive. Moreover,
from (\ref{010908b}) and Theorem~\ref{260808a} we obtain the following statement.

\thrml{010908c}
Let $\A$ be an S-ring over a CG-ring.
Then the ring $\A$ is a non-trivial generalized wreath product if and only if so
is the ring $\wh\A$. More exactly, $\A$ satisfies the $I/J$-condition if and only
if $\wh\A$ satisfies the $J^\bot/I^\bot$-condition.\bull
\ethrm

The following theorem shows that an S-ring and its dual are cyclotomic
or not simultaneously. This fact can also be deduced from the results of~\cite{GC92} by
using the well-known 1-1 correspondence between S-rings and translation
association schemes.

\thrml{091208a}
Let $\A=\cyc(K,R)$ where $K\le R^\times$. Then $\wh\A=\cyc(\wh K,\wh R)$.
\ethrm
\proof Let $X\in\orb(\wh K,\wh R)$. Then given $\chi_1,\chi_2\in X$ there
exists $r\in K$ such that $\chi^{}_1=\chi_2^{\,(r)}$. Since $S=rS$ for
each basic set $S$ of~$\A$, this implies that
$$
\chi^{}_1(S)=\chi_2^{\,(r)}(S)=\chi^{}_2(rS)=\chi^{}_2(S),\qquad S\in\S(\A).
$$
Therefore $\A'\ge\wh\A$ where $\A'=\cyc(\wh K,\wh R)$. On the other hand,
$\rk(\wh\A)=\rk(\A)$ (see Subsection~\ref{151209u}). Moreover, since there is a
ring isomorphism from $R$ to $\wh R$ taking $K$ to $\wh K$ (see Subsection~\ref{021109k}),
we also have $\rk(\A)=\rk(\A')$. Thus $\rk(\wh\A)=\rk(\A')$ whence it follows that
$\wh\A=\A'$.\bull

\section{Rational basic set outside a maximal $\A$-ideal}\label{241109a}

In this section we prove Theorem~\ref{231109a} from which Theorem~\ref{201109a} immediately
follows. Below a subset $X$ of a ring $R$ is called {\it $Q$-rational} for some
$Q\subset\P(R)$ if $R_Q^\times X=X$. When $Q=\{p\}$, we say that $X$ is $p$-rational.

\thrml{231109a}
Let $\A$ be an S-ring over a CG-ring  $R$ and $I\in\I_{max}(\A)$. Suppose
that any basic set in $\S(\A)_{R\setminus I}$ is $Q$-rational where $Q=\P(R/I)$.
Set $I^\star$ to be the intersection of all $\A$-ideals containing~$R_Q$.
Then the ring~$\A$ satisfies the $I/J$-condition where $J=I\cap I^\star$. Moreover, if $J=0$,
then $\A=\A_{I^{}}\otimes \A_{I^\star}$.
\ethrm

The proof of the Theorem~\ref{231109a} will be given in the end of the section. In the following
auxiliary statement we study some properties of the $\A$-ideals $I$, $I^\star$ and $J$
defined in it. We keep the notations and the hypothesis of this
theorem except for the $Q$-rationality.

\lmml{011009d}
The following statements hold:
\nmrt
\tm{1} $I+I^\star=R$, $I^\star=R_Q+J_{Q'}$ and $|R/I|=|I^\star/J|=|R_Q/J_Q|$,
\tm{2} $\I_{max}(\A_{I^\star})=\{J\}$,
\tm{3} for each $p\in Q$ the ideal $I^\star$ equals the intersection of all $\A$-ideals
containing~$R_p$,
\tm{4} $\LI(X)_Q\subset J$ for all $X\in\S(\A)_{R\setminus I}$.
\enmrt
\elmm
\proof By the definition of $Q$ we have $R_Q\not\subset I$. Therefore
$I^\star\not\subset I$, and $I+I^\star=R$ by the maximality of~$I$. This gives
the natural group isomorphism
\qtnl{271009a}
I^\star/J\cong R/I
\eqtn
and hence $|R/I|=|I^\star/J|$. Since $\P(R/I)=Q$, it follows that  $I^\star=R_Q+J_{Q'}$
and $|I^\star/J|=|R_Q/J_Q|$ which completes the proof of statement~(1).
Isomorphism (\ref{271009a}) induces the bijection
$$
\{I'\in\I(R):\ J\subset I'\subset I^\star\}\to
\{I'\in\I(R):\ I\subset I'\subset R\},\quad
I'\mapsto I+I'
$$
the converse to which takes $I'$ to
$I'\cap I^\star$. This induces a bijection on the $\A$-ideals. Thus statement~(2)
follows from the maximality of~$I$. To prove statement~(3) let $p\in Q$.
Denote by $I'$ the intersection of all $\A$-ideals containing $R_p$. Then
obviously $I'\subset I^\star$. By statement~(2) this implies that either
$I'\subset J$ or $I'=I^\star$. Since $p\in\P(I^\star/J)$ (see statement~(1)) and
$R_p\subset I'$, the inclusion is impossible and we are done. Let us prove statement~(4).
Since $X\in\S(\A)$ and $I\in\I_{max}(\A)$, we see that $\LI(X)\subset I$.
On the other hand, $\LI(X)_Q\subset R_Q\subset I^\star$. Thus
$\LI(X)_Q\subset I\cap I^\star=J$.\bull\medskip

In fact, Theorem~\ref{231109a} will be deduced from the following statement where we
use the same notations.

\thrml{011009a}
Let $\A$ be an S-ring over a CG-ring $R$, $I\in\I_{max}(\A)$ and
$X\in\S(\A)_{R\setminus I}$. Suppose that the set $X$ is $p$-rational
for some $p\in Q$. Then $J\subset\LI(X)$ and
\qtnl{061009c}
\pi(X)=\pi(X)_{Q}+\pi(X)_{Q'}
\eqtn
where $\pi=\pi_J$. Moreover, both $\pi(X)_Q$ and $\pi(X)_{Q'}$ are basic sets of
the S-ring $\A_{R/J}$, and $\pi(X)_Q=\pi(I^\star)^\#$.
\ethrm
\proof By statement~(4) of Lemma~\ref{011009d} without loss of generality we
can assume that $\LI(X)_p=0$. First, we claim that
\qtnl{231109b}
\rk(\A_{R/I})=2.
\eqtn
Indeed, the maximality of $I$ implies that the S-ring $\A_{R/I}$ is $R/I$-primitive.
Therefore~(\ref{231109b}) follows from Theorem~\ref{201207a} whenever the ring~$R/I$
is not a field. However, if it is a field, then $\P(R/I)=\{p\}$ and~(\ref{231109b})
follows from the $p$-rationality of~$X$.\medskip

Let us prove that $J\subset\LI(X)$. To do this we observe that from~(\ref{231109b})
it follows that $\pi_I(X)=\pi_I(R_Q)^\#$. So $X_p\ne \{0\}$. Since $R_p^\times X=X$, this
implies by statement~(3) of Theorem~\ref{261009a} that $\LI(X\cup Y)_p\ne 0$
where $Y$ is as in this theorem. Therefore
\qtnl{011009c}
\LI(X\cup Y)\supset I^\star_0
\eqtn
where $I^\star_0$ is the intersection of all $\A$-ideals $I'$ with $I'_p\ne 0$.
Clearly, $I^\star_0\in\I(\A)$. Moreover, from the definition of~$I^\star$ it follows that
$I^\star_0\subset I^\star$. Therefore by statement~(2) of Lemma~\ref{011009d}
we have
\qtnl{021009d}
I^\star_0=I^\star\quad\text{or}\quad I^\star_0\subset J.
\eqtn
Let $I'\subset I^\star_0$ be a non-zero $\A$-ideal of $R$ other than $I^\star$.
Then from (\ref{011009c}) it follows that $Y+I'\subset X\cup Y$.
On the other hand,
\qtnl{021009a}
X\cap (Y+I')=\emptyset.
\eqtn
Indeed, otherwise $X\subset Y+I'$ because $X$ is a basic set of~$\A$
and $Y+I'$ is an $\A$-set. Since $X\in\S(\A)_{R\setminus I}$ and
$\rk(\A_{R/I})=2$, there exists $x\in X$ such that $x_p\in R_p^\times$.
However, $(X^{[p]})_p=\{0\}$ by statement~(1) of Theorem~\ref{261009a}.
So $Y_p=\{0\}$. Moreover, $I'\subset J$ by statement~(2) of
Lemma~\ref{011009d}. Therefore from statement~(1) of this lemma it follows
that $R_p^\times\cap I'=\emptyset$. Thus $x\not\in I'$. Contradiction.
Thus (\ref{021009a}) holds, and hence $Y+I'=Y$. But then
using (\ref{011009c}) we obtain that $I'\subset\LI(X)$. This implies that
\qtnl{021009e}
I_0^\star=I^\star.
\eqtn
Otherwise, $I_0^\star\subset J\ne I^\star$ and we can take $I'=I^\star_0$,
which contradicts the assumption that $\LI(X)_p=0$. In this case we can take
$I'=J$, which proves that $J\subset\LI(X)$.\medskip

To complete the proof it suffices to assume that $J=0$. Then $I^\star=R_Q$ and
$I=R_{Q'}$. Since both of them are $\A$-ideals, we have $X_Q,X_{Q'}\in\S(\A)$
(Lemma~\ref{130209d}). Moreover, from statement~(2) of Lemma~\ref{011009d} it follows that
$X_Q=(R_Q)^\#$. Next, by (\ref{011009c}) and (\ref{021009e}) the set
$X\cup Y$ is a union of $R_Q$-cosets. On the other hand, by the definition
of $Y$ we have $Y_{p'}\subset X_{p'}$, and hence
$Y_{Q'}\subset X_{Q'}$. Thus
\qtnl{021002e}
X\cup Y=R_Q+X_{Q'}.
\eqtn
However, since $X_Q=(R_Q)^\#$, we have $X\subset R_Q^\#+X_{Q'}$. On the
other hand, $Y\subset X_{Q'}$, for otherwise $Y_Q\ne \{0\}$, and
hence $Y_Q\supset R_Q^\#$ which contradicts the fact that
$Y_p=\{0\}$. Thus from (\ref{021002e}) it follows that $X=R_Q^\#+X_{Q'}$.\bull\medskip

{\bf Proof of Theorem~\ref{231109a}.} Since the set $\S(\A)_{R\setminus I}$
consists of $Q$-rational sets, we conclude by Theorem~\ref{011009a} that
$J\subset\LI(X)$ for all $X\in\S(\A)_{R\setminus I}$. Thus the S-ring~$\A$
satisfies the $I/J$-condition. To complete the proof suppose that $J=0$.
Then $I^\star=R_Q$ and $I=R_{Q'}$ where $Q'=\P(R)\setminus Q$. Therefore from
Theorem~\ref{011009a} it follows that for all $X\in\S(\A)_{R\setminus I}$
we have $X_{Q^{}}\in\S(\A_{I})$, $X_{Q'}\in\S(\A_{I^\star})$ and $X=X_{Q^{}}+X_{Q'}$.
Since these three relations obviously hold  for $X\in\S(\A)_I$, we are done.\bull

\section{S-rings with $\I_{max}(\A)\ne\I_{max}(R)$}\label{021209u}

It is known that any non-dense S-ring over a cyclic group is either a non-trivial generalized
wreath product, or a tensor product one factor of which is an S-ring of rank~$2$
(see~\cite[Theorem~5.3]{EP01be}). In the case of S-rings over an arbitrary CG-ring we can prove
the same statement only under a stronger condition yet.

\thrml{140409f}
Let $\A$ be an S-ring over a CG-ring $R$. Suppose that $\I_{max}(\A)\ne\I_{max}(R)$.
Then $\A$ is either a non-trivial generalized wreath product, or a tensor product one factor of
which is an S-ring of rank~$2$ over a non-field.
\ethrm

Since $I\in\I_{max}(\A)$ is a maximal ideal of $R$ if and only if $R/I$ is a field,
Theorem~\ref{140409f} is an immediate consequence of the following statement in which
we keep the notations of Section~\ref{241109a}.

\thrml{231009a}
Let $\A$ be an S-ring over a CG-ring  $R$ and $I\in\I_{max}(\A)$. Suppose that
$R/I$ is not a field. Then the ring~$\A$ satisfies the $I/J$-condition. Moreover,
if $J=0$, then $\A=\A_{I^{}}\otimes \A_{I^\star}$.
\ethrm

In its turn Theorem~\ref{231009a} is an immediate consequence of Theorem~\ref{231109a}
and the following statement which will be proved a bit later.

\thrml{061009a}
Let $\A$ be an S-ring over a CG-ring $R$ and $I\in\I_{max}(\A)$.
Suppose that $R/I$ is not a field. Then any basic set in $\S(\A)_{R\setminus I}$
is $Q$-rational where $Q=\P(R/I)$.
\ethrm

\proof We need a special consequence of Theorem~\ref{011009a} that gives us a convenient
form of a $Q$-rational basic set.

\lmml{061009b}
In the conditions of Theorem~\ref{011009a} we have $X=(R_Q\setminus J_Q)+X_{Q'}.$
\elmm
\proof Since the full $\pi$-preimages of $\pi(X)$ and $\pi(X)_Q$ are $X$ and $I^\star\setminus J$
respectively, from the equality (\ref{061009c}) we obtain that
\qtnl{061009d}
X=(I^\star\setminus J)+Y
\eqtn
where $Y$ is the full preimage of $\pi(X)_{Q'}$. Taking into account that
$Q\cap Q'=\emptyset$, we have $\pi(Y_Q)=\pi(Y)_Q=0$, and hence $Y_Q\subset J_Q$.
Therefore
$$
J_Q+y=J_Q+y_Q+y_{Q'}=J_Q+y_{Q'}
$$
for all $y\in Y$. Since
$Y$ is a union of $J_Q$-cosets, this shows that $Y=J_Q+Y_{Q'}$.
On the other hand, $X_{Q'}$ and $Y_{Q'}$ are unions of $J_{Q'}$-cosets,
and so $J_Q+X_{Q'}$ and $J_Q+Y_{Q'}$ are unions of $J$-cosets. Since also
$\pi(X_{Q'})=\pi(Y_{Q'})$, we have $J_Q+X_{Q'}=J_Q+Y_{Q'}$. Thus from
(\ref{061009d}) we obtain
\qtnl{061009e}
X=(I^\star\setminus J)+Y=(I^\star\setminus J)+J_Q+Y_{Q'}=(I^\star\setminus J)+J_Q+X_{Q'}=(I^\star\setminus J)+X_{Q'}.
\eqtn
Finally, from statement~(1) of Lemma~\ref{011009d} it follows that
$$
I^\star\setminus J=(R_Q+J_{Q'})\setminus (J_Q+J_{Q'})=(R_Q\setminus J_Q)+J_{Q'}.
$$
However, $X_{Q'}$ is a union of $J_{Q'}$-cosets. Thus we are done
due to (\ref{061009e}).\bull\medskip

Turn to the proof of Theorem~\ref{061009a}. Let $X\in\S(\A)_{R\setminus I}$ and $p\in Q$.
Set $Y=R^\times_pX$. Then obviously $Y_{Q'}=X_{Q'}$. Therefore by Lemma~\ref{061009b}
we have
\qtnl{061009f}
Y=(R_Q\setminus J_Q)+X_{Q'}.
\eqtn
If $|Q|>1$, then by statement~(5) of Lemma~\ref{011009d} there exists
an element of $R_Q\setminus J_Q$ with zero $p$-coordinate. Due to (\ref{061009f})
one can find an element in $Y$, and hence in $X$, say $x$, with the same
property. Then $x_p=0$ and hence $R_p^\times x=x$. Therefore $R_p^\times X=X$.
Thus the basic set $X$ is $Q$-rational and we are done. In the remaining case
$|Q|=1$ we make use of the following lemma proved in Section~\ref{241109c}.

\lmml{091009a}
Suppose that $Q=\{p\}$. Then
\nmrt
\tm{1} if $p\not\in\P(J)$, then $R_p^\times X=X$,
\tm{2} if $p\in\P(J)$, then $\LI(X)\cap J\ne 0$.
\enmrt
\elmm

To complete the proof let $Q=\{p\}$. Then by statement~(1) of Lemma~\ref{091009a}
we can assume that $p\in\P(J)$. In this case by statement~(2) of that lemma
$J_0=\LI(X)\cap J$ is a non-zero $\A$-ideal of~$R$. Set $R'=\pi(R)$,
$\A'=\A_{R'}$, $I'=\pi(I)$ and $X'=\pi(X)$ where $\pi=\pi_{J_0}$. Then
obviously $I'\in\I(\A')$ and $X'\in\S(\A')_{R'\setminus I'}$. Moreover, since
$J_0\subset J\subset I$ we have
$$
\P(R'/I')=\P(R/I)=Q=\{p\}.
$$
Finally, $\LI(X')\cap J'=0$, and hence $p\not\in\P(J')$ by statement~(2) of
Lemma~\ref{091009a}. Thus from statement~(1) of that lemma it follows that
$(R'_p)^\times X'=X'$. Since $X$ is a union of $J_0$-cosets, this implies that
$R_p^\times X=X$.\bull

\section{Decomposition of a pure S-ring}\label{021209i}

We note that not every pure S-ring over a CG-ring is dense, e.g. take the S-ring of
rank~$2$ over a non-field. The following theorem shows
that at least in the odd case this is essentially  a unique reason for
a pure S-ring not to be dense.

\thrml{160409a}
Any pure S-ring over a CG-ring of odd characteristic
is the tensor product of a dense pure S-ring and S-rings of rank~$2$ over
non-fields.
\ethrm

We will prove Theorem~\ref{160409a} in the end of this section. In what follows
we say that a pure S-ring is {\it indecomposable} if it is not a tensor product, one factor
of which is an S-ring of rank~$2$ over a non-field.

\thrml{301009a}
Let $\A$ be an S-ring over a CG-ring $R$. Then the following statements are
equivalent:
\nmrt
\tm{1} $\A$ is pure indecomposable,
\tm{2} $\wh\A$ is pure indecomposable,
\tm{3} $\A$ is pure and $\I_{max}(\A)=\I_{max}(R)$,
\tm{4} $\A$ is pure and $\I_{min}(\A)=\I_{min}(R)$.
\enmrt
\ethrm
\proof From Theorem~\ref{010908c} it follows that if $\A$ or $\wh\A$ is pure,
then neither $\A$ nor $\wh\A$ is a nontrivial generalized wreath product.
Besides, it is easily seen that if $\A$ is a tensor product one factor of
which is an S-ring of rank~$2$ over a non-field, then $\I(\A)\ne\I(R)$.
By Theorem~\ref{140409f} these facts prove the equivalence
$(1)\,\Leftrightarrow\,(3)$. To complete the proof of the theorem it suffices to
verify the implication $(1)\,\Rightarrow\,(2)$. Indeed, if it is true, then the
converse implication follows by duality, and the equivalence
$(3)\,\Leftrightarrow\,(4)$ is an immediate consequence of the fact that
$\I_{min}(R)=\I_{min}(\A)$ if and only if $\I_{max}(\wh R)=\I_{max}(\wh\A)$
(see~(\ref{010908b})).\medskip

To prove the implication $(1)\,\Rightarrow\,(2)$ suppose that $\A$ is a pure
indecomposable S-ring. Then by statement~(2) of Theorem~\ref{140509a} and
Theorem~\ref{140409f} it suffices to verify that $\wh\A$ is a pure S-ring. For
this purpose we note that due to the implication $(1)\,\Rightarrow\,(3)$, the set
$R^\times$ being the complement in $R$
to the union of all maximal $\A$-ideals, is an $\A$-subset of $R$. This implies
that the basic set $X$ of $\A$ that contains $1$ is a subset of $R^\times$. Thus by
Theorem~\ref{230810b} the set $X$ is a subgroup of $R^\times$, which is pure by the hypothesis. Let
$\chi\in\wh R^\times$ and $\wh X$ the basic set of $\wh\A$ containing $\chi$. It
suffices to verify that $\LI(\wh X)=0$. Suppose that this is not true. Then
\qtnl{060209b}
\wh X\chi'=\wh X,\qquad \chi'\in\LI(\wh X).
\eqtn
Since the set $X$ is pure, by statement (2) of Theorem~\ref{250303b}
for $S=X$ we have $\chi(X')\ne 0$ where $X'=Xr$ for some
$r\in R^\times$. Since $X'\in\S(\A)$ from (\ref{060209b}) we obtain that
$\chi(X')=\chi\chi'(X')$ for all $\chi'\in\LI(\wh X)$. So
\qtnl{060209c}
a=\sum_{\chi'\in\LI(\wh X)}\chi\chi'(X')\ne 0.
\eqtn
On the other hand, set $K=1+I$ where $I$ is the preimage of the ideal
$\LI(\wh X)$ with respect to the isomorphism from $R$ to $\wh R$ induced by~$\chi$.
Then when $\chi'$ runs over the set $\LI(\wh X)$ the element $\chi\chi'$ runs over
the set $\{\chi^{(r)}:\ r\in K\}$. Then
$$
a=\sum_{r\in K}\chi^{(r)}(X')=\sum_{r\in K}\sum_{x\in X'}\chi^{(r)}(x)
=\sum_{x\in X'}\sum_{r\in K}\chi(rx).
$$
However $\chi(I)=0$ because $\chi$ is non-trivial on~$I$. Therefore for all
$x\in R^\times$ we have
$$
\sum_{r\in K}\chi(rx)=\sum_{h\in I}\chi(x+hx)=
\chi(x)\sum_{h\in I}\chi(hx)=
\chi(x)\chi(I)=0.
$$
Thus $a=0$ because $X'\subset R^\times$,  which contradicts~(\ref{060209c}).\bull\medskip

\thrml{060209a}
The  S-ring dual to a pure S-ring over a CG-ring is also pure.
\ethrm
\proof By statement~(2) of Theorem~\ref{140509a} without loss of generality
we can assume that the input S-ring is indecomposable. Then the required
statement immediately follows from the implication $(1)\,\Rightarrow\,(2)$
of Theorem~\ref{301009a}.\bull\medskip

{\bf Proof of Theorem~\ref{160409a}.} By statement~(2) of Theorem~\ref{140509a} without loss of generality
we can assume that the S-ring~$\A$ is indecomposable. Therefore Theorem~\ref{160409a}
is an immediate consequence of the following statement.

\thrm
Suppose that the characteristic of a CG-ring~$R$ is odd. Then any pure
indecomposable S-ring over~$R$ is dense.
\ethrm
\proof Let $\A$ be a pure indecomposable S-ring over the ring~$R$.
Suppose that $R_p$ is not a
field for some $p$. Denote by $J$ the minimal ideal of $R_p$. Then by the
implication $(1)\,\Rightarrow\,(4)$ of Theorem~\ref{301009a} we have
$J\in\I(\A)$. Since obviously $J\ne R_p$ and the characteristic of
the ring~$R$ is odd, from Theorem~\ref{021109a} it follows that the S-ring $\A_{R/J}$
is pure. Moreover, it is easily seen that any maximal ideal of $R/J$
is the $\pi_J$-image of a maximal ideal of~$R$. Therefore by the
equivalence $(1)\,\Leftrightarrow\,(3)$ of Theorem~\ref{301009a} the S-ring~$\A_{R/J}$
is a pure indecomposable one. So by induction we conclude that
$\I(\A_{R/J})=\I(R/J)$. Thus $\I(\A)$ contains any ideal $I$ such that
$I_p\ne 0$ for some $p\in\P$ for which $R_p$ is not a field. Therefore by the
implication $(1)\,\Rightarrow\,(4)$ of Theorem~\ref{301009a} we have
$\I(\A)=\I(R)$.\bull

\section{Dense pure S-rings}\label{271109a}

In this section we prove the following theorem which provides together with
Theorem~\ref{160409a} a straightforward deduction of Theorem~\ref{270409a}.

\thrml{130209b}
Any dense pure S-ring over a CG-ring of odd characteristic is cyclotomic.
\ethrm

The proof will be given in the end of the section throughout which we fix
a dense S-ring $\A$ over a CG-ring $R$. From equality~(\ref{010908b}) it
follows that the ring $\wh\A$ is also dense. We need two
lemmas in each of which the characteristic of the ring~$R$ is arbitrary.

\lmml{270308a}
Suppose that $\orb(K,R^\times)\subset\S^*(\A)$ where  $K\le R^\times$. Then any
pure orbit of the group $\wh K$ in $\wh R$ belongs to~$\S^*(\wh\A)$.
\elmm
\proof Let $X_1$ be a pure orbit of the group $\wh K$. Then $X_1=\chi_1^{\wh K}$
for some character $\chi_1\in\wh R$. Denote by $X$ the basic set of
$\wh\A$ containing $\chi_1$ and set $X_2=\chi_2^{\wh K}$ where
$\chi_2\in X$. Since $\wh\A$ is dense, the set $Y=\wh R^\times\chi_1$
belongs to $\S^*(\wh\A)$ (see~(\ref{021109b})). This
implies that $X\subset Y$ whence it follows that
$\chi_1,\chi_2\in Y$, and hence $X_1,X_2\subset Y$.\medskip

Denote by $a$ the cardinality of the kernel of the natural action of the group $K$
on the set $Y$. Since the action is semiregular (see Section~\ref{021109d}),
given $S\in\orb(K,R^\times)$ and $s\in S$ we have
$$
\chi_i(S)=\sum_{r\in K}\chi_i(rs)=\sum_{r\in K}\chi_i^{(r)}(s)=as(X_i),\quad i=1,2,
$$
where $s(X_i)$ is defined by (\ref{201108b}) with $G=\wh R^+$, $S=X_i$ and
$\chi$ being the character of $G$ corresponding to~$s$.
On the other hand, as $S\in\S^*(\A)$ the definition of the dual S-ring implies
that $\chi_1(S)=\chi_2(S)$. Thus
$$
s(X_1)=s(X_2),\quad s\in R^\times.
$$
Since $X_1$ and $X_2$ are pure orbits of the group $K$, from statement~(1) of
Theorem~\ref{250303b} applied to $\wh R$ and $X_1,X_2$ it follows that
$X_1=X_2$ and hence $\chi_2\in X_1$. However, $\chi_2$ is an arbitrary
element of~$X$. Thus $X\subset X_1$ and we are done.\bull\medskip

From now on we assume that for any $p\in\P$ where $\P=\P(R)$, the characteristic of the
Galois ring $R_p$ equals $p^{n_p}$ for some $n_p\ge 1$. Set
\qtnl{261109d}
\P'=\{p\in\P:\ n_p>1\},\qquad R'=\bigcup_{\P'\subset Q\subset\P}R^\times_Q.
\eqtn
Due to the density of~$\A$, from (\ref{021109b}) it follows that
$R^\times_Q$ is $\A$-set for all~$Q$, and hence $R'\in\S^*(\A)$. Moreover, by
Lemma~\ref{130209d} and the definition of~$\P'$ we have
\qtnl{261109a}
\S(\A)_{R'}=\{X_Q:\ X\in\S(\A)_{R^\times},\ Q\supset\P'\}.
\eqtn
Below for $p\in\P'$ we denote by $m(p)$ the product of all $q\in\P'\setminus\{p\}$.

\lmml{120209a}
Let $\S(\A)_{R^\times}=\orb(K,R^\times)$ where $K\le R^\times$. Suppose
that for some $p\in\P'$ the orbits of the group~$K$ on the set $pR^\times\cup m(p)R^\times$
are pure $\A$-sets. Then $\S(\A)_{pR^\times}=\orb(K,pR^\times)$.
\elmm
\proof Suppose on the contrary that there is an orbit $X\in\orb(K,pR^\times)$
containing two distinct basic sets $Y$ and $Y'$ of~$\A$. Then
\qtnl{261109c}
Y'=rY
\eqtn
for some $r\in K$. By the hypothesis of the lemma the sets $X$ and $Z=mK$ with $m=m(p)$,
are pure. Therefore, by Theorem~\ref{250303b} (with $S=Y$ and $S'=Y'$ for
statement (1), and with $S=Z$ for statement~(2)) there exist a character
$\chi\in\wh R$ and a set $Z'=r'Z$ with $r'\in R^\times$ such that
\qtnl{170303c}
\chi(Y)\ne \chi(Y'),\quad \chi(Z')\ne 0.
\eqtn
However, by the definition of $m$ we have $Y+Z'\subset R'$ where $R'$ is the
$\A$-set defined in~(\ref{261109d}). Since also
$Y\in\S(\A)$ and $Z'\in\orb(K,mR^\times)\subset\S^*(\A)$, we conclude that
\qtnl{130209a}
\xi(Y)\xi(Z')\in\A_{R'}.
\eqtn
Taking into account that $\S(\A)_{R^\times}=\orb(K,R^\times)$ and that
the set of $K$-orbits is closed with respect to taking projections, we see
by (\ref{261109a}) that the right-hand side of~(\ref{130209a}) is $K$-invariant.
This implies that so is the left-hand side. Therefore by (\ref{261109c}) and the
definition of $Z'$ we have
$$
\xi(Y)\xi(Z')=r(\xi(Y)\xi(Z'))=\xi(rY)\xi(rZ')=\xi(Y')\xi(Z').
$$
Applying $\chi$ to both sides of this equality we obtain a contradiction
with~(\ref{170303c}).\bull\medskip

{\bf Proof of Theorem~\ref{130209b}.} Let $\A$ be a dense pure S-ring over
a CG-ring $R$ of odd characteristic. Due to the density of~$\A$, Theorem~\ref{230810b}
implies that there exists a group $K\le R^\times$ such that
\qtnl{091208d}
\S(\A)=\orb(K,R^\times).
\eqtn
So by Lemma~\ref{270308a} any pure orbit of the group~$\wh K$
belongs to~$\S^*(\wh\A)$. However, since the S-ring $\A$ is pure, the group $K$
and hence the group~$\wh K$ are also pure. By Theorem~\ref{090209a}
this implies that
\qtnl{170209a}
\orb(\wh K,m\wh R)\subset \S^*(\wh\A)_{m\wh R}
\eqtn
for all integers $m$ dividing $n/\prod_{p\in\P}p$. However, by Lemma~\ref{130209d}
each basic set of $\wh\A$ is the projection of some element of $\S(\wh\A)_{m\wh R}$
with $m$ as above. Therefore inclusion (\ref{170209a}) holds for all~$m$.
Thus  $\wh\A\ge\cyc(\wh K,\wh R)$, and hence $\A\ge\cyc(K,R)$ by
Theorem~\ref{091208a}. Since the group $K$ is pure,
this implies by~(\ref{091208d}) and  Lemma~\ref{120209a} that
$$
\S(\A)_{pR^\times}=\orb(K,pR^\times),\qquad p\in\P',
$$
where $\P'$ is defined as in~(\ref{261109d}). Moreover, by Theorem~\ref{090209a}
the group $K_{pR}=f_{pR}(K)$ where $f_{pR}$ is the epimorphism defined in the end of
Subsection~\ref{021109dd} (see (\ref{161109a}) and below),
is a pure subgroup of the group $(R_{pR})^\times$ (here $K_{pR}=pK$ as sets,
see~(\ref{261109u})). So by induction we can assume that
\qtnl{261109r}
\S(\A)_{pR}=\cyc(K_{pR},R_{pR}),\qquad p\in\P'.
\eqtn
However, $\S(\A)=\S(\A)_{R'}\cup\S'$ where $R'$ is as in~(\ref{261109d}) and
$\S'=\bigcup_{p\in\P'}\S(\A)_{pR}$. Besides, due to~(\ref{091208d})
and~(\ref{261109a}) we have $\S(\A)_{R'}\subset\orb(K,R)$. Since by (\ref{261109r})
we also have $\S'\subset\orb(K,R)$, it follows that $\S(\A)\subset\orb(K,R)$.
Thus $\A=\cyc(K,R)$.\bull

\section{Proof of Theorem~\ref{210809a}}\label{191109a}

Denote by $\T_p$ and $\U_p$ (resp. by $\T_q$ and $\U_q$) the Teichm\"uller
group and the group of principal units of the ring~$R_p$ (resp.~$R_q$).
Since $q$ divides $p^d-1$, the cyclic group $\T_p$ of order $p^d-1$ contains
a unique subgroup~$T_p$ of order~$q$. Similarly, since $p$ divides $q^e-1$,
the cyclic group $\T_q$ of order $q^e-1$ contains a unique subgroup~$T_q$ of
order~$p$. Let us consider a subgroup of $R^\times$ defined as follows:
$$
K=(T_p\U_p)\times(T_q\U_q).
$$
Thus $K_p=T_p\U_p$ and $K_q=T_q\U_q$ (see Notation). Any group $L\le K$ with $L_p=K_p$ and
$L_q=K_q$ (i.e. a subdirect product of~$K_p$ and~$K_q$) is determined by means
of an appropriate group $L_0$ and epimorphisms $f_p:K_p\to L_0$, $f_q:K_q\to L_0$ as follows:
$$
L=\{(u,v)\in K:\ f_p(u)=f_q(v)\}.
$$

Let us define subgroups $K_1$ and $K_2$ of the group $K$ in the above way where $f_p$ and $f_q$ are fixed
epimorphisms on a cyclic group $L_0$ of order $q$ in the first case and
of order $p$ in the second one. Clearly, the sets $R^\times$ and $pR\cup qR$
are both $K_1$-invariant and $K_2$-invariant. Set
\qtnl{210709a}
\S=\orb(K_1,R^\times)\cup\orb(K_2,pR\cup qR).
\eqtn
Then to prove Theorem~\ref{210809a} it suffices to verify the following statement.

\thrml{260509d}
The $\Z$-module $\A=\Span_\Z\{\xi(X):\ X\in\S\}$ is a non-pure dense S-ring over
the ring $R$. Moreover, this S-ring cannot be a non-trivial generalized wreath
product.
\ethrm

Before giving the proof of Theorem~\ref{260509d} let us cite some simple
properties of the groups~$K$, $K_1$ and~$K_2$. The following statement is
straightforward.

\lmml{260509a}
There exist direct decompositions $\U_p=U_p\cdot U'_p$ and $\U_q=U_q\cdot U'_q$
with cyclic groups $U_p$ and $U_q$, such that
\qtnl{301109a}
K_1=(\U_p\times T_q U'_q)\cdot L_1, \qquad
K_2=(T_p U'_p\times\U_q)\cdot L_2
\eqtn
where $L_1$ is a subdirect product of $T_p$ and $U_q$ of order $q$
and $L_2$ is a subdirect product of $U_p$ and $T_q$ of order $p$,
and both decompositions (\ref{301109a}) are direct.\bull
\elmm

From Lemma~\ref{260509a} it follows that
\qtnl{011209a}
|K_1|=p^{d+1}q^e,\ \LI(K_1)=pR_p,\qquad
|K_2|=p^dq^{e+1},\ \LI(K_2)=qR_q.
\eqtn
Moreover,
$$
K_1\cdot K_2=K,\qquad K_1\cap K_2=(U'_p\times U'_q)\cdot L_1\cdot L_2
$$
where the latter decomposition is direct, and for $i,j\in\{0,1,2\}$ we have
\qtnl{250509a}
\orb(K,p^iq^jR^\times)=\orb(K_1,p^iq^jR^\times)=\orb(K_2,p^iq^jR^\times),\qquad i+j\in\{2,3,4\},
\eqtn
whereas if $(i,j)=(0,1)$ or $(i,j)=(1,0)$, then
\qtnl{250509b}
\orb(K,qR^\times)=\orb(K_1,qR^\times),\qquad \orb(K,pR^\times)=\orb(K_2,pR^\times).
\eqtn

{\bf Proof of Theorem~\ref{260509d}.} Obviously, the elements of the set $\S$
form a partition of~$R$ such that any ideal of~$R$ is a union of classes
of the partition. Moreover, the induced partition of~$R^\times$ consists of
orbits of the group~$K_1$ which is not pure. Since also $R^\times\S=\S$,
to prove the first of the theorem it suffices to verify that $\A$ is an S-ring
over the group $R^+$.\medskip

From the definition of~$\S$ it follows that $\{0\}\in\S$ and $-\S=\S$. Therefore we
need to check only that $\xi(X)\xi(Y)\in\A$ for all $X,Y\in\S$. We distinguish
three cases depending on to which parts of~$\S$ the sets $X$ and $Y$ belong. Denote
by $\xi_1$ and $\xi_2$ the elements of $\Z R$ such that
$$
\xi(X)\xi(Y)=\xi_1+\xi_2,\quad \supp(\xi_1)\subset R^\times,\quad
\supp(\xi_2)\subset pR\cup qR.
$$
Below for $\xi\in\Z R$ we set $\LI(\xi)$ to be the largest ideal $I\in\I(R)$ for which
$\xi\xi(I)=|I|\xi$; obviously, if $\xi=\xi(Z)$ for some $Z\subset R$, then
$\LI(\xi)=\LI(Z)$.\medskip

{\bf Case 1: $X,Y\in\orb(K_1,R^\times)$}. In this case $\xi(X),\xi(Y)\in\cyc(K_1,R)$,
and hence $\xi_1,\xi_2\in\cyc(K_1,R)$. Thus since $K\ge K_2$, it suffices to verify that
$\xi_2\in\cyc(K,R)$. However, since obviously
$\xi(qR)\in\cyc(K_1,R)$, from~(\ref{250509a}) and~(\ref{250509b}) it follows that
$\xi_2\circ\xi(qR)\in\cyc(K,R)$. It remains to show that the element
$\xi=\xi_2\circ\xi(pR)$ belongs to $\cyc(K,R)$. To do this we observe that by~(\ref{011209a})
we have $\LI(X)=\LI(Y)=pR_p$. Therefore $\LI(\xi(X)\xi(Y))\ge pR_p$
and hence
$$
\LI(\xi)\ge pR_p.
$$
This implies that $\xi=\xi(pR_p)\xi'$ for some $\xi'\in\Z R_q$. In particular, all
elements of each $pR_p$-coset enter $\xi$ with the same coefficient. Therefore taking
into account that $\xi\in\cyc(K_1,R)$, we conclude that
$\xi'\in \cyc(K,R)$. Since also $\xi(pR_p)\in\cyc(K,R)$, we obtain that
$\xi\in\cyc(K,R)$.\medskip

{\bf Case 2: $X,Y\in\orb(K_2,pR\cup qR)$}. Arguing as above we obtain that
$\xi_1,\xi_2\in\cyc(K_2,R)$. Therefore without loss of generality we can assume that
$\xi_1\ne 0$ (see~(\ref{210709a})). Then it is easily seen that $\xi_2=0$.
Moreover, we can assume that $X\subset p^iR^\times$
and $Y\subset q^jR^\times$ where $i,j\in\{1,2\}$. It suffices to verify that
\qtnl{210709b}
\xi_X,\xi_Y\in\cyc(K,R),\qquad \LI(\xi_Y)=pR_p
\eqtn
where $\xi_X=\xi(pR_p)\xi(X)$ and $\xi_Y=\xi(qR_q)\xi(Y)$. Indeed, since
$\LI(K_2)=qR_q$, it follows that $\LI(X)\ge qR_q$. The converse inclusion
is obvious. Thus $\LI(X)=qR_q$. Therefore
using the second part of (\ref{210709b}) we obtain that
$$
\xi(X)\xi(Y)=(c_q\xi(X)\xi(qR_q))\xi(Y)=
c_q\xi(X)(\xi(qR_q)\xi(Y))=c_q\xi(X)\xi_Y=
$$
$$
c_q\xi(X)(c_p\xi(pR_p)\xi_Y)=c_pc_q(\xi(pR_p)\xi(X))\xi_Y=c_pc_q\xi_X\xi_Y
$$
where $c_p=|pR_p|^{-1}$ and $c_q=|qR_q|^{-1}$. Thus $\xi_X\xi_Y\in\cyc(K,R)$ by
the first part of (\ref{210709b}).\medskip

To prove (\ref{210709b}) let $X=p^iuK_2$ for some $u\in R^\times$. By Lemma~\ref{260509a}
we have $K_2=(T_pU'_p\times \U_q)\cdot L_2$ and $L_2=\{(f(t),t):\ t\in T_q\}$ where
$f:T_q\mapsto U_p$ is a group isomorphism. Therefore
$$
p^iL_2=\{(p^if(t),p^it):\ t\in T_q\}=\{p^i\}\times p^iT_q,
$$
and
$$
X=(u_pT_pU'_p\times u_q\U_q)\cdot p^iL_2=
(u_pT_pU'_p\times u_q\U_q)\cdot(\{p^i\}\times p^iT_q)=
p^iu_pT_p\times p^iu_qT_q\U_q.
$$
So $\xi_X=\xi(pR_p)\xi(X)=c\,\xi(pR_p\times p^iu_qT_q\U_q)$ for a positive integer $c$.
Thus $\xi_X\in\cyc(K,R)$.\medskip

Next, let $Y=q^juK_2$ where
$u\in R^\times$. Then as above we have
\qtnl{040809a}
Y=(q^ju_pT_pU'_p\times q^ju_q\U_q)\cdot L_2=\bigcup_{t\in T_q}Z_t\times\{q^ju_qt\}
\eqtn
where $Z_t=q^ju_pf(t)T_pU'_p$. Since obviously
$\xi(qR_q)\xi(Z_t\times\{q^ju_qt\})=\xi(Z_t\times qR_q)$, this implies that
$$
\xi_Y=\xi(qR_q)\xi(Y)=\sum_{t\in T_q}\xi(Z_t\times qR_q)=
$$
$$
\sum_{r\in U_p}\xi(q^ju_prT_pU'_p\times qR_q)=\xi(q^ju_pT_p\U_p\times qR_q).
$$
Thus $\xi_Y\in\cyc(K,R)$. Since $q^ju_pT_p\U_p=q^ju_pT_p+pR_p$, we also obtain that
$\LI(\xi_Y)=pR_p$.\medskip

{\bf Case 3: $X\in\orb(K_1,R^\times)$, $Y\in\orb(K_2,pR\cup qR)$} (or vice versa).
First, suppose that $Y\in\orb(K_2,pR)$. Then from~(\ref{250509a}) and the second equality of~(\ref{250509b}) it follows
that $Y\in\orb(K,pR)$, and hence $\xi(Y)\in\cyc(K_1,R)$. This implies
that $\xi(X)\xi(Y)\in\cyc(K_1,R)$. Therefore $\xi_1\in\A$. Since also
$\supp(\xi_2)\subset qR$, by~(\ref{250509a}) and the first equality of~(\ref{250509b})
we conclude that $\xi_2\in\cyc(K,R)\subset\A$. Thus $\xi(X)\xi(Y)\in\A$.\medskip

Suppose that $Y\in\orb(K_2,qR\setminus pR)$. Then $Y=q^juK_2$ where
$u\in R^\times$ and $j\in\{1,2\}$. The element $Z_t$ defined in Case~2 can be
rewritten in the following form: $Z_t=\bigcup_{s\in T_p}(r_ts+r_tsH)$ where
$r_t=qu_pf(t)$ and $H=U'_p-1\subset pR_p$. However, for any $s\in T_p$ we have
$$
\xi(r_ts+r_tsH)\xi(pR_p)=|H|\xi(r_ts+pR_p)=|H|\xi(r_ts\U_p)=|H|\xi(qu_ps\U_p).
$$
Together with (\ref{040809a}) this implies that
$$
\xi(Y)\xi(pR_p)=
\xi(\bigcup_{t\in T_q}Z_t\times\{qu_qt\})\xi(pR_p)=
$$
$$
(\sum_{t\in T_q}\xi(Z_t\times\{qu_qt\})\xi(pR_p)=
\sum_{t\in T_q}\sum_{s\in T_p}\xi(r_ts+r_tsH)\xi(pR_p)\xi(qu_qt)=
$$
$$
|H|\sum_{s\in T_p}\sum_{t\in T_q}\xi(qu_ps\U_p)\xi(qu_qt)=
|H|\xi(qu_pT_p\U_p\times qu_qT_q).
$$
Thus $\xi(Y)\xi(pR_p)\in\cyc(K,R)$. Besides, since $\LI(X)=pR_p$, we have
$$
\xi(X)\xi(Y)=c\,\xi(X)(\xi(Y)\xi(pR_p))
$$
where $c$ is a positive rational. Therefore $\xi(X)\xi(Y)\in\cyc(K_1,R)$,
whence it follows that $\xi_1,\xi_2\in\cyc(K_1,R)$. In particular,
$\xi_1\in\A$. To prove that $\xi_2\in\A$, it suffices to verify that $\xi_2\in\cyc(K,R)$.
However, it is easy to see that $\supp(\xi_2)\subset pR$. Moreover, $\LI(\xi_2)\ge pR_p$
because $\supp(\xi_1)\cap pR=\emptyset$ and $\LI(\xi(X)\xi(Y))\ge\LI(X)=pR_p$.
Thus applying to $\xi=\xi_2$ the same argument as in the end of Case~2, we obtain that
$\xi_2\in\cyc(K,R)$ which completes the proof of Case~3 and the first part of
the theorem.\medskip

To prove the second part of the theorem suppose on the contrary that the S-ring~$\A$
satisfies the $I/J$-condition non-trivially. Then without loss of generality
we can assume that $I$ is a maximal $\A$-ideal of~$R$, i.e. $I=pR$ or $I=qR$.
Since $R^\times\subset R\setminus I$ and $\LI(X)=pR_p$ for all $X\in\S(\A)_{R^\times}$
due to~(\ref{011209a}), we also can assume that $J=pR_p$. Suppose that $I=pR$.
Then the set $X=qK_2$ belongs to $\S(\A)_{R\setminus I}$, and hence $J\le \LI(X)$. So
$\LI(X)_p\ne 0$. On the other hand, a straightforward computation based on
Lemma~\ref{260509a} shows that $\LI(X)_p=0$. Contradiction. Finally,
suppose that $I=qR$. Then the set $X=pK_2$ belongs to $\S(\A)_{R\setminus I}$, and hence
$\LI(X)_p\ne 0$. On the other hand,
$\LI(X)_p=0$ because $X$ cannot contain any $pR_p$-coset. Contradiction.\bull

\section{Proof of Theorem~\ref{250303b}}\label{271109b}

The proof is based on the following two results on abelian groups to be proved
in the end of the section. Let us fix some notation for an abelian group~$G$.
Denote by $G_0$ the socle of~$G$; in our case it is the product of all subgroups of~$G$
of prime order.
For $p\in\P$ and $Q\subset\P$ where $\P=\P(G)$, we set $G_{0,p}=(G_0)_p$ and $G_{0,Q}=(G_0)_Q$.

\thrml{200409b}
Let $G$ be an abelian group of exponent~$m$ each Sylow $p$-subgroup of which is
homocyclic and let $K\subset\C$ be a field linearly separated from $\Q[w]$ over $\Q$ where $w$
is a primitive $m$th root of unity. Denote by $\Psi=\Psi_K(G)$ the set of all
$K$-epimorphisms $\psi:KG\to K[w]$. Then
\qtnl{020209a}
\bigcap_{\psi\in\Psi}\ker(\psi)=I_K(G)
\eqtn
where $I_K(G)=\sum_{p\in\P}I_p\otimes KG_{p'}$ and $I_p$ is the ideal of $KG_p$ spanned
by $\xi(G_{0,p})$.\footnote{Here and below all tensor products are taken over~$K$.}
\ethrm

\lmml{050209a}
In the notation of Theorem~\ref{200409b} given a nonzero $\xi\in I_K(G)$
there exists a nonempty set $Q\subset\P$ and a $G_{0,Q}$-coset $A\subset G$
such that $\supp(\xi_A)_p\supset A_p$ for all $p\in Q$ where
$\xi_A=\xi\circ\xi(A)$.
\elmm

Let us turn to the proof of Theorem~\ref{250303b}. By the hypothesis
there exist a group $L\le R^\times$ and a pure orbit $T$ of it
such that $S,S'\subset T$.
To prove statement (1) suppose on the contrary that $\chi(S)=\chi(S')$ for
all $\chi\in\wh R^\times$. However, when a character $\chi$ runs over $\wh R^\times$ its
extension $\psi:\Q G\to\C$ where $G=R^+$, runs over the set $\Psi_\Q(G)$ defined
in Theorem~\ref{200409b}. Thus $\xi(S)-\xi(S')\in\ker(\psi)$ for all $\psi\in\Psi_\Q(G)$.
Since each Sylow subgroup of~$G$ is homocyclic in our case, by Theorem~\ref{200409b} this implies that
$\xi(S)-\xi(S')\in I_K(G)$. Let us prove that this contradicts the
purity of~$T$.\medskip

Since $S\ne S'$ the element $\xi=\xi(S)-\xi(S')$ is nonzero.
So by Lemma~\ref{050209a} there exist
a nonempty set $Q\subset\P$ and a $G_{0,Q}$-coset $A\subset G$  such
that $\supp(\xi_A)_p\supset A_p$ for all $p\in Q$. Since $S,S'\subset T$, this
implies that
\qtnl{210409a}
(T\cap A)_p=A_p,\qquad p\in Q.
\eqtn
One can see that $T\cap A$ is a block of the abelian group $L$ acting on $T$. This
implies that $T\cap A\in\orb(L')$
where $L'$ is the setwise stabilizer of $T\cap A$ in~$L$. Therefore
\qtnl{210409b}
|T\cap A|=|L_A|
\eqtn
where $L_A$ is the permutation group on $T\cap A$ induced by~$L'$. It is easily
seen that given $p\in Q$ the family $M_p$ of all nonempty sets $X\cap T$,
$X\in A/G_{0,Q\setminus\{p\}}$, forms an imprimitivity system
for $L_A$. From (\ref{210409a}) it follows that $|M_p|=|G_{0,p}|$ for all $p\in Q$.
Besides, due to (\ref{210409b}) the number $|M_p|$ divides $|T\cap A|$
for all $p\in Q$.  Thus
$$
|T\cap A|\ge\prod_{p\in Q}|G_{0,p}|=|A|
$$
whence it follows that $T\supset A$. Therefore $T=T+G_{0,Q}$, and hence the set $T$ is not pure. This
contradiction completes the proof of statement~(1).\medskip

To prove statement (2) let $\chi\in\wh R^\times$. By statement~(1) with
$S'=\emptyset$ there exists a character~$\chi'\in\wh R^\times$ such that
$\chi'(S)\ne\chi'(S')=0$. However, due to~(\ref{301008b}) we have
$\chi'=\chi^{(r)}$ for some $r\in R^\times$. Thus, $\chi(rS)=\chi'(S)\ne 0$.\bull\medskip

{\bf Proof of Theorem~\ref{200409b}.} Clearly, the right-hand side of (\ref{020209a}) is contained in the
left-hand side. Let us prove the converse inclusion by induction on $|\P|$.
Without loss of generality we assume that $|\P|>0$. Fix $q\in\P$ and set $G'=G_{q'}$
Then each $g\in G$ can uniquely be written in the form $g=g'g_q$ where $g'\in G'$ and $g_q\in G_q$.\medskip

For induction purposes some preliminary work is needed. Set $K'=K[w_q]$ where $w_q$ is a primitive $m_q$th
root of unity. For $\xi=\sum_{g\in G}a_gg$ belonging to~$KG$ and $\psi_q\in\Psi_K(G_q)$
define an element $\xi'=\xi'(\psi_q)$ of the ring $K'G'$ by
$$
\xi'=\sum_{g'\in G'}a'_{g'}g'\quad\text{with}\quad
a'_{g'}=\sum_{g_q\in G_q}a_{g'g_q}\psi_q(g_q).
$$
Next, using the equalities $KG=KG'\otimes KG_q$ and $KG'\otimes K'=K'G'$ let us define a ring
epimorphism
$$
\wt\psi_q:KG\to K'G'
$$
by $\wt\psi_q=\id_{KG'}\otimes\psi_q$. Then it is easily seen that
\qtnl{251209c}
\wt\psi_q(\xi)=\xi',\quad \ker(\wt\psi_q)=KG'\otimes\ker(\psi_q),\quad
\wt\psi_q(I_K(G')\otimes KG_q)=I_{K'}(G').
\eqtn

To prove that the left-hand side of~(\ref{020209a}) is contained in the right-hand side
suppose that $\xi\in \bigcap_{\psi\in\Psi}\ker(\psi)$. First, we claim that
\qtnl{251209a}
\xi'(\psi_q)\in I_{K'}(G')\quad\text{for all}\ \psi_q\in\Psi_K(G_q).
\eqtn
Indeed, let $\psi_q\in\Psi_K(G_q)$. We note that
the field $K'$ is linearly separated from $\Q[w']$ over $\Q$ where $w'$ is a
primitive $m_{q'}$th root of unity because $K$ is linearly separated from $\Q[w]$ over $\Q$.
Therefore by the induction hypothesis for the group $G'$ and the field $K'$
it suffices to check that $\psi'(\xi')=0$ for all $\psi'\in\Psi'$ where  $\xi'=\xi'(\psi_q)$
and $\Psi'=\Psi_{K'}(G')$. However,
$$
\psi'(\xi')=\psi'(\sum_{g'\in G'}a'_{g'}g')=
\psi'(\sum_{g'\in G'}\sum_{g_q\in G_q}a_{g'g_q}\psi_q(g_q)g')=
$$
$$
\sum_{g'\in G'}\sum_{g_q\in G_q}a_{g'g_q}\psi'(g')\psi_q(g_q)=
\sum_{g'\in G'}\sum_{g_q\in G_q}a_{g'g_q}\psi(g'g_q)=
\sum_{g\in G}a_g\psi(g)=\psi(\xi)
$$
where $\psi=\psi'_{KG'}\otimes\psi_q\in\Psi$ with $\psi'_{KG'}$ being the restriction of
$\psi'$ to $KG'$). By the choice of $\xi$, we conclude that $\psi'(\xi')=\psi(\xi)=0$,
and we are done.\medskip

Further, set $M=I_K(G')\otimes KG_q$ and $N_{\psi_q}=KG'\otimes\ker(\psi_q)$ where
$\psi_q\in\Psi_K(G_q)$. Then from (\ref{251209c}) and (\ref{251209a}) it immediately follows that
for all $\psi_q$ we have
$$
\xi\in(\wt\psi_q)^{-1}(\xi')\subset
(\wt\psi_q)^{-1}(I_{K'}(G'))\subset
M+\ker(\wt\psi_q)=M+N_{\psi_q}.
$$
Therefore to complete the proof, i.e. to verify that $\xi\in I_K(G)$ it suffices to show that
\qtnl{030209a}
\bigcap_{\psi_q\in\Psi_K(G_q)}(M+N_{\psi_q})=I_K(G).
\eqtn

To prove (\ref{030209a}) we observe that
$$
(M+N_{\psi_q})/M=
N_{\psi_q}/(M\cap N_{\psi_q})=
$$
$$
N_{\psi_q}/(I_K(G')\otimes\ker(\psi_q))=(KG'/I_K(G'))\otimes\ker(\psi_q).
$$
On the other hand, from \cite{EP09} it follows that
$\bigcap_{\psi_q\in\Psi_K(G_q)}\ker(\psi_q)=I_K(G_q)$. Thus,
\qtnl{030209b}
\bigcap_{\psi_q\in\Psi_K(G_q)}(M+N_{\psi_q})/M=
(KG'/I_K(G'))\otimes I_K(G_q).
\eqtn
Similarly,
\qtnl{030209c}
I_K(G)/M=(M+KG'\otimes I_K(G_q))/M=(KG'/I_K(G'))\otimes I_K(G_q).
\eqtn
Therefore (\ref{030209a}) follows from (\ref{030209b}) and (\ref{030209c}).\bull\medskip

{\bf Proof of Lemma~\ref{050209a}.}
Denote by $M$ the set of all $Q\subset\P$ for which there exists
a $G_{0,Q}$-coset $A$ such that
\qtnl{271109g}
\xi_A\in\Span_K\{\xi(a+G_{0,p}):\ a\in A,\ p\in Q\}^\#.
\eqtn
From the definition of $I_K(G)$ it follows that
$\xi$ is a $K$-linear combination of elements $\xi(a+G_{0,p})$ with $a\in G$ and
$p\in\P$. Since also $\xi\ne 0$, we see that $\P\in M$. Let $Q$ be a minimal
(by inclusion) subset of $\P$ belonging to $M$ and $A$ the corresponding coset.
Then obviously $Q\ne\emptyset$. We claim that for $Q$ and $A$ the statement
of the theorem holds. Suppose on the contrary that this is not true. Then
there exists $p\in Q$ such that $\supp(\xi_A)_p\supsetneq A_p $. This implies that
\qtnl{200409a}
\supp(\xi_A)\cap B=\emptyset
\eqtn
for some coset $B\in A/G_{0,Q\setminus\{p\}}$. However, it is easily
seen that given $b\in B$ and $q\in Q\setminus\{p\}$ we have
$$
\xi(b+G_{0,q})=
\sum_{x\in G_{0,q}}\xi(b+x+G_{0,p})-
\sum_{y\in G_{0,p}^\#}\xi(b+y+G_{0,q}).
$$
Therefore, without loss of generality we can assume that in a representation
of $\xi_A$ afforded by~(\ref{271109g}) the coefficient at $\xi(b+G_{0,q})$ equals $0$
for all $b\in B$. Due to (\ref{200409a}) this implies that the coefficient at
$\xi(a+G_{0,p})$ equals $0$ for all $a\in A$. Thus condition (\ref{271109g})
holds for $Q$ and $A$ replaced by $Q\setminus\{p\}$ and any $G_{0,p}$-coset in $A$ that
intersects $\supp(\xi_A)$, respectively. But this contradicts the
minimality of~$Q$.\bull\medskip

\section{Proof of Lemma~\ref{091009a}}\label{241109c}

From the assumptions it follows that $R_p=\GR(p^n,d)$ where $n\ge 2$.

\lmml{161009a}
Suppose that $R_q\subset J$ for some $q\in\P(R)$. Then $R_q\subset\LI(X)$.
\elmm
\proof Clearly, $q\ne p$. By Theorem~\ref{011009a} the set $R_p^\times X$ is a union
of~$J$-cosets. Due to the assumption each of them contains an element with zero
$q$-coordinate. So the set $X$ also contains such an element. Therefore
$R_q^\times X=X$, and hence by statement~(2) of Theorem~\ref{261009a}
(for $p=q$)
$$
\LI(X)_q\ne 0\quad\text{or}\quad X^{[q]}\ne\emptyset.
$$
In the former case $(\LI(X)\cap J)_q=\LI(X)_q\ne 0$ and we are done by induction.
Let us show that the latter case is impossible. Indeed, let us consider
the $\A$-ideal
$$
I'=\UI(X^{[q]})\cap I^\star.
$$
Since $X\subset R\setminus I$ and $Q=\{p\}$, we have $x_p\not\in J$ for
all $x\in X$. However, $q\ne p$ and hence $x_p\not\in J$ for all $x\in X^{[q]}$.
This implies that $\UI(X^{[q]})_p\not\subset J_p$. Since $I^\star\supset R_p$, it
follows that $I'_p\not\subset J_p$. Therefore $I'\not\subset J$. By
statement~(2) of Lemma~\ref{011009d} this shows that $I'=I^\star$. Thus
$I'_q=I^\star_q\supset J_q=R_q$. On the other hand, taking into account that
$(X^{[q]})_q=\{0\}$ (statement~(1) of Theorem~\ref{261009a}), we see that $I'_q=0$.
Contradiction.\bull\medskip

\crllrl{271009d}
If $p\not\in\P(J)$, then $I^\star\setminus J\in\S(\A)$.
\ecrllr
\proof It is easily seen that the hypothesis of Lemma~\ref{091009a} is satisfied
for $\A$, $R$ and $I$ replaced with $\A_{I^\star}$, $R_{I^\star}$ and $J$
(see Lemma~\ref{011009d}). Moreover, since
$p\not\in\P(J)$, we have $(R_{I^\star})_q\subset J$ for all $q\in\P(R_{I^\star})$,
$q\ne p$. By Lemma~\ref{161009a} this implies that $(R_{I^\star})_q\subset\LI(X)$ for
all~$q$ and any $X\in\S(\A_{I^\star\setminus J})$. In particular,
$J\le\LI(X)$. Since by statement~(2) of Lemma~\ref{011009d} and Theorem~\ref{201207a}
we have $\rk(\A_{I^\star/J})=2$, this implies that $X=I^\star\setminus J$
and we are done.\bull\medskip

Let us turn to the proof of Lemma~\ref{091009a}. To prove statement~(1) suppose that
$p\not\in\P(J)$. Denote by $t$ the cardinality of the set
$\{uX:\ u\in R_p^\times\}$. Then it suffices to verify that $t=1$. To do this
we observe that from Lemma~\ref{061009b} it follows that for any $x'\in X_{p'}$ the
cardinality of the set
$$
F_{x'}=\{x_p:\ x\in X,\ x_{p'}=x'\}
$$
equals $f=|R_p^\#|/t$. Next, by Corollary~\ref{271009d} we have
$I^\star\setminus J\in\S(\A)$, and hence
$$
\xi(X)\xi(-X)=\alpha\xi(I^\star\setminus J)+\xi'
$$
where $\alpha$ is a non-negative integer and $\xi'\circ\xi(I^\star\setminus J)=0$. It
is easily seen that the element $x-y$ with
$x,y\in X$ belongs to $I^\star\setminus J$ if and only if $x_p\ne y_p$ and
$x_{p'}-y_{p'}\in J$. Since $|I^\star\setminus J|=|J||R_p^\#|$, it follows that
\qtnl{140409c}
\alpha\le f^2|X_{p'}||J|/|I^\star\setminus J|=f|X_{p'}|/t.
\eqtn
On the other hand, taking into account that $p\not\in\P(J)$ one can find
$x_0\in X$ such that $(x_0)_p$ is of order $p$. Since
$(1+\rad(R_p))\{x_0\}=\{x_0\}$, we have $(1+\rad(R_p))X=X$. Besides, by
Lemma~\ref{061009b} it is also true that
\qtnl{161009v}
|F_{x'}\cap p^iR_p|=|p^iR_p^\times|/t,\qquad i=0,\ldots,n-1,\ x'\in X_{p'}.
\eqtn
Therefore the set $F_{x'}\setminus p^{n-1}R_p$
is a union of $p^{n-1}R_p$-cosets. This implies that given $x'\in X_{p'}$ an
element of $I^\star\setminus J$ of order $p$ can be represented as $x-y$ with $x,y\in X$,
$x_{p'}=x'$ in at least $|F_{x'}\setminus p^{n-1}R_p|=f-f_0$ ways where
$f_0=|F_{x'}\cap p^{n-1}R_p|=|(p^{n-1}R_p)^\#|/t$ (see~(\ref{161009v})). This
implies that
$$
\alpha\ge |X_{p'}|(f-f_0).
$$
Together with (\ref{140409c}) this shows that $f\ge tf-tf_0$ whence it follows
that $|R_p^\#|\ge t(|R_p|-|p^{n-1}R_p|)$.
Since $|R_p|=p^{nd}$ and $|p^{n-1}R_p|=p^d$, this gives
$p^{nd}>t(p^{nd}-p^d)$. For $t\ge 2$ we have $p^{(n-1)d}<t/(t-1)\le 2$,
which is impossible for $n\ge 2$. Thus $t=1$ and statement~(1) is proved.\medskip

To prove statement~(2) suppose that $p\in\P(J)$. By Lemma~\ref{161009a} we can assume that
$R_q\not\subset J$ for all $q\in\P(R)$. It suffices to verify that
\qtnl{151009a}
(1+\rad(R_p))X=X.
\eqtn
Indeed, we have $I_p=J_p$ by the lemma hypothesis and the definition of~$J$, and
$J_p\ne 0$ because $p\in\P(J)$.
Therefore $I_p\ne 0$, and hence $X_p\subset R_p\setminus I_{0,p}$.
Due to (\ref{151009a}) this implies that the set $X$ is a union
of $I_{0,p}$-cosets. So $\LI(X)_p\ne 0$ and we are done because
$J_p\supset I_{0,p}$. To prove (\ref{151009a}) we note that by
the assumption $1+J$ is a subgroup of~$R^\times$. Then the set
$$
Y=(1+J)X
$$
belongs to $\S^*(\A)$. It is easily seen that the hypothesis of Lemma~\ref{091009a}
is satisfied for $\A$ and $X$ replaced with $\A_{R/J}$ and $\pi_J(Y)$. So by already
proved statement~(1) of this lemma we have
$\pi_J(Y)=\pi_J(Y')$ for all $Y'\in M$ where $M=\{rY:\ r\in R_p^\times\}$.
Therefore for any $y\in R$ the
sets $Y^{}_{J,y}$ and $Y'_{J,y}$ are empty or not simultaneously.
By Lemma~\ref{090608a} this implies that
$$
|Y^{}_{J,y}|=|Y'_{J,rY}|=|Y'_{J,y}|
$$
where $Y'=rY$. On the other hand, by Theorem~\ref{011009a} applied to
$R_p^\times Y$ the full $\pi_J$-preimage of $\pi_J(R^\times_p Y)$ coincides with
the union of all sets from~$M$. This implies that given $y\in Y$ the set
$y+J$ is a disjoint union of the sets $Y'_{J,y}$, $Y'\in M$. Thus
\qtnl{151009c}
|J|=|y+J|=\sum_{Y'\in M}|Y'_{J,y}|=|M||Y_{J,y}|.
\eqtn
Next, take $y\in Y$ such that $y_p\in R^\times_p$ (such an $y$ exists
because $X\in\S(\A)_{R\setminus I}$ and $Q=\{p\}$). Since
$J_p\subset\rad(R_p)$, it follows that $y'_p\in R^\times_p$ for all $y'\in Y_{J,y}$.
This implies that $|(1+J_p)y'|=|J_p|$ for all~$y'\in Y_{J,y}$. Taking into account that
$Y_{J,y}$ is a union of the sets $(1+J_p)y'$, we conclude that $|J_p|$ divides
$|Y_{J,y}|$. By~(\ref{151009c}) this means that $|M|$ is coprime to~$p$.
Since the abelian group $R_p^\times$ acts transitively on~$M$, this
implies that $1+\rad(R_p)$ is a subgroup of the kernel of this action.
Therefore $(1+\rad(R_p))Y=Y$, and hence
$$
(1+\rad(R_p))(1+J)X=(1+\rad(R_p))Y=Y=(1+J)X.
$$
In particular, given $u\in 1+\rad(R_p)$ and $z\in X$ one can find
$u'\in 1+J$ and $z'\in X$ such that $uz=u'z'$. It follows that $(u/u')X=X$.
When $u$ runs over the elements of highest order in $1+\rad(R_p)$,
the element $(u/u')_p$ runs a full system of representatives of $(1+J_p)$-cosets
in $1+\rad(R_p)$. Therefore the Sylow $p$-subgroup of the group $K$ generated
by all elements $u/u'$ coincides with $1+\rad(R_p)$. Since $KX=X$, the equality
(\ref{151009a}) holds and we are done.

\end{document}